\documentclass[a4paper,11pt]{article}
\usepackage{latexsym,amssymb}
\usepackage{amsmath}
\usepackage[mathscr]{eucal}
\usepackage{color}
\usepackage[abbrev]{amsrefs}

\usepackage{latexsym}
\usepackage{amsmath}
\usepackage{amssymb}
\usepackage{enumerate}
\usepackage{mathrsfs}
\usepackage{color}

\newtheorem{Theorem}{\bf Theorem}[section]
\newtheorem{Lemma}{\bf Lemma}[section]
\newtheorem{Proposition}{\bf Proposition}[section]
\newtheorem{Corollary}{\bf Corollary}[section]
\newtheorem{Remark}{\bf Remark}[section]
\newtheorem{Example}{\bf Example}[section]
\newtheorem{Definition}{\bf Definition}[section]

\newenvironment{theorem}{\begin{Theorem}$\!\!\!$}{\end{Theorem}}
\newenvironment{lemma}{\begin{Lemma}$\!\!\!$}{\end{Lemma}}
\newenvironment{proposition}{\begin{Proposition}$\!\!\!$}{\end{Proposition}}
\newenvironment{corollary}{\begin{Corollary}$\!\!\!$}{\end{Corollary}}
\newenvironment{remark}{\begin{Remark}$\!\!\!$}{\end{Remark}}

\newenvironment{definition}{\begin{Definition}$\!\!\!$}{\end{Definition}}

\numberwithin{equation}{section}

\def\XXint#1#2#3{{\setbox0=\hbox{$#1{#2#3}{\int}$}
\vcenter{\hbox{$#2#3$}}\kern-.5\wd0}}

\begin{document}

\title{Initial traces and solvability \\for a semilinear heat equation
on a half space of ${\mathbb R}^N$}
\author{
\qquad\\
Kotaro Hisa, Kazuhiro Ishige, and Jin Takahashi
}
\date{}
\maketitle

\begin{abstract}
We show the existence and the uniqueness of initial traces of nonnegative solutions to a semilinear heat equation on 
a half space of ${\mathbb R}^N$ under the zero Dirichlet boundary condition. 
Furthermore, we obtain necessary conditions and sufficient conditions on the initial data 
for the solvability of the corresponding Cauchy--Dirichlet problem. 
Our necessary conditions and sufficient conditions are sharp 
and enable us to find optimal singularities of initial data for the solvability of the Cauchy--Dirichlet problem. 
\end{abstract}
\vspace{25pt}
\noindent Addresses:

\smallskip
\noindent K.~H.:  Mathematical Institute, Tohoku University,\\
\qquad\,\,\,\,\, 6-3 Aoba, Aramaki, Aoba-ku, Sendai 980-8578, Japan. \\
\noindent 
E-mail: {\tt kotaro.hisa.d5@tohoku.ac.jp}\\

\smallskip
\noindent 
K.~I.:  Graduate School of Mathematical Sciences, The University of Tokyo,\\
\qquad\,\,\, 3-8-1 Komaba, Meguro-ku, Tokyo 153-8914, Japan. \\
\noindent 
E-mail: {\tt ishige@ms.u-tokyo.ac.jp}\\

\smallskip
\noindent J.~T.:  Department of Mathematical and Computing Science,\\ 
\qquad\,\,\, Tokyo Institute of Technology, 2-12-1 Ookayama, Meguro-ku,\\ 
\qquad\,\,\, Tokyo 152-8552, Japan.\\
\noindent 
E-mail: {\tt takahashi@c.titech.ac.jp}\\
\vspace{20pt}

\noindent
{\it MSC:} 35K58, 
35A01, 35A21, 35K20 
\vspace{3pt}

\noindent
{\it Keywords:} initial trace, semilinear heat equation, Cauchy--Dirichlet problem, 
solvability
\vspace{3pt}
\date{}
\maketitle

\newpage
\section{Introduction}
Let $T\in(0,\infty]$, $\Omega:={\mathbb R}^{N-1}\times(0,\infty)$ if $N\ge 2$, and $\Omega:=(0,\infty)$ if $N=1$. 
This paper is concerned with initial traces of nonnegative, measurable, and finite almost everywhere functions 
in $\Omega\times(0,T)$ solving the problem
\begin{equation}
\label{eq:E}
\tag{E}
\left\{
\begin{array}{ll}
\partial_t u=\Delta u+u^p,\quad & x\in\Omega,\,\,\,t\in(0,T),\vspace{3pt}\\
u(x,t)=0, & x\in\partial \Omega,\,\,\,t\in(0,T),
\end{array}
\right.
\end{equation}
where $p>1$, 
and the solvability of the corresponding Cauchy--Dirichlet problem. 
Generally, 
qualitative properties of initial traces of nonnegative solutions for parabolic equations have been studied 
in the framework of nonnegative Radon measures and they have played important roles 
in the study of necessary conditions on the solvability for the corresponding initial value problems. 
On the other hand, in our problem~(E), 
due to the zero Dirichlet boundary condition, 
initial traces of nonnegative solutions 
cannot be treated in the framework of nonnegative Radon measures on $\overline{\Omega}$
and little is known concerning qualitative properties of initial traces. 
Indeed, in the case of $1<p<1+2/(N+1)$, 
there exists a positive smooth function~$v$ solving problem~(E) for some $T>0$ with 
\begin{equation}
\label{eq:1.1}
v(x,0)=-\partial_{x_N}\delta_N\quad\mbox{on}\quad\overline{\Omega},
\end{equation}
where $\delta_N$ is the $N$-dimensional Dirac measure concentrated at the origin. (See e.g. \cites{KY, TW}. See also Corollary~\ref{Corollary:5.3}.) 
Then the initial condition of the solution~$v$ is not a nonnegative Radon measure on $\overline{\Omega}$.

In this paper, in order to overcome the difficulty due to the zero Dirichlet boundary condition, 
for any nonnegative, measurable, and finite almost everywhere function~$u$ solving problem~(E), 
we propose to consider the initial trace of the function $x_Nu(x,t)$, 
instead of the function~$u$ itself. 
We prove the existence and the uniqueness of the initial trace of the function $x_Nu(x,t)$ 
in nonnegative Radon measures on $\overline{\Omega}$. 
Furthermore, we obtain necessary conditions on the initial data 
for the solvability of the corresponding Cauchy--Dirichlet problem
 \begin{equation}
  \tag{P}
  \label{eq:P}
  \left\{
  \begin{array}{ll}
  \partial_t u=\Delta u+u^p,\quad & x\in\Omega,\,\,\,t>0,\vspace{3pt}\\
  u=0, & x\in\partial \Omega,\,\,\,t>0,\vspace{3pt}\\
  x_Nu(x,0)=\mu, & x\in\overline{\Omega},
  \end{array}
  \right.
  \end{equation}
  where $\mu$ is a nonnegative Radon measure on $\overline{\Omega}$.
More precisely, 
\begin{itemize}
  \item[(O1)] 
  for any nonnegative, measurable, and finite almost everywhere function~$u$ in $\Omega\times(0,T)$ solving problem~(E), 
  we show the existence and the uniqueness of 
  a nonnegative Radon measure $\nu$ on ${\mathbb R}^N$ 
  with $\mbox{supp}\,\nu\subset\overline{\Omega}$ satisfying
  $$
  \underset{t\to +0}{\mbox{ess lim}}\int_\Omega x_N u(x,t)\phi(x)\,dx=\int_{\overline{\Omega}}\phi(x)\,d\nu(x)
  \quad\mbox{for all $\phi\in C_0(\overline{\Omega})$};
  $$
  \item[(O2)] 
  we formulate a definition of the solution to Cauchy--Dirichlet problem~\eqref{eq:P} in $\Omega\times(0,T)$, 
  and we show that the initial data $\mu$ coincides with the initial trace of $x_Nu(x,t)$;
  \item[(O3)] 
  we obtain necessary conditions on the initial data
  for the existence of solutions to problem~\eqref{eq:P}.
\end{itemize}
In our setting, initial condition~\eqref{eq:1.1} is regarded as 
$x_Nu(x,0)=\delta_N$ on $\overline{\Omega}$. 
(See also Remark~\ref{Remark:1.1} and Corollary~\ref{Corollary:5.3}.) 
Moreover, 
\begin{itemize}
  \item[(O4)] 
  we obtain sufficient conditions for the solvability of problem~\eqref{eq:P}, 
  and find optimal singularities of the initial data for the solvability of problem~\eqref{eq:P}. 
\end{itemize}  
Then we see that it is valid to consider the initial trace of $x_Nu(x,t)$, instead of $u$, 
and optimal singularities given in (O4) show that our necessary conditions and sufficient conditions for solvability are sharp. 
In our arguments, the explicit representation of the Dirichlet heat kernel in $\Omega\times(0,\infty)$ is crucial.
\vspace{3pt}

The study of initial traces of nonnegative solutions to parabolic equations is a classical subject
and it has been investigated for various parabolic equations,
for example, the heat equation (see \cites{A, W}), the porous medium equation (see \cites{AC, BCP, HP}),
the parabolic $p$-Laplace equation (see \cites{DH, DH02}), the doubly nonlinear parabolic equation (see \cites{I, IJK, ZX}),
the fractional diffusion equation (see \cite{BSV}), the Finsler heat equation (see \cite{AIS}), 
and parabolic equations with nonlinear terms (see e.g.~\cites{ADi, BP, BCV, FHIL, FI01, HI18, HI19, IKO, MV, MV02, MV03, Zhao}). 
Among others, 
in \cite{HI18} the first author and the second author of this paper proved  
the existence and the uniqueness of 
initial traces of nonnegative solutions to the semilinear heat equation
\begin{equation}
\label{eq:E'}
\tag{E'}
\partial_t u=\Delta u+u^p,\quad  x\in{\mathbb R}^N,\,\,\,t\in(0,T),
\end{equation}
where $p>1$. Furthermore, they studied necessary conditions and sufficient conditions 
for the solvability of the Cauchy problem 
\begin{equation}
\label{eq:P'}
\tag{P'}
\left\{
\begin{array}{ll}
\partial_t u=\Delta u+u^p,\quad & x\in{\mathbb R}^N,\,\,\,t>0,\vspace{3pt}\\
u(x,0)=\eta, & x\in{\mathbb R}^N,
\end{array}
\right.
\end{equation}
where $\eta$ is a nonnegative Radon measure on ${\mathbb R}^N$. 
These necessary conditions and sufficient conditions are sharp and they depend on whether $1<p<p_N$, $p=p_N$, or $p>p_N$. 
Here and in what follows, $p_d:=1+2/d$ for $d=1,2,\dots$. 
\subsection{Problems~\eqref{eq:E'} and \eqref{eq:P'}}
Let $T\in(0,\infty)$. 
Let $u$ be a nonnegative, measurable, and finite almost everywhere function 
in ${\mathbb R}^N\times(0,T)$ solving problem~\eqref{eq:E'}, 
that is, 
$u$ satisfies 
$$
u(x,t)=\int_{{\mathbb R}^N}\Gamma_N(x-y,t-\tau)u(y,\tau)\,dy
+\int_\tau^t\int_{{\mathbb R}^N}\Gamma_N(x-y,t-s)u(y,s)^p\,dy\,ds
$$
for almost all (a.a.)~$(x,t)\in{\mathbb R}^N\times(0,T)$ and a.a.~$\tau\in(0,T)$ with $\tau<t$. 
Here $\Gamma_N$ is the heat kernel in ${\mathbb R}^N\times(0,\infty)$ (see \eqref{eq:1.3}). 
It follows from \cite{HI18}*{Theorem~1.1} (see also \cite{FHIL}*{Corollary~3.1}) that 
there exists a unique nonnegative Radon measure $\theta$ on ${\mathbb R}^N$ satisfying
$$
\underset{t\to +0}{\mbox{ess lim}}\int_{{\mathbb R}^N}u(x,t)\phi(x)\,dx=\int_{{\mathbb R}^N}\phi(x)\,d\theta(x)
$$
for all $\phi\in C_0({\mathbb R}^N)$. 
Furthermore, 
\begin{itemize}
  \item[(F1)] 
  there exists $\gamma_1=\gamma_1(N,p)>0$ such that 
  \begin{equation}
  \label{eq:1.2}
  \sup_{x\in{\mathbb R}^N}\theta(B(x,\sigma))\le\gamma_1\sigma^{N-\frac{2}{p-1}}
  \quad
  \mbox{for all $\sigma\in(0,\sqrt{T})$}.
  \end{equation}
\end{itemize}
Here $B(x,\sigma):=\{y\in{\mathbb R}^N:|x-y|<\sigma\}$ for $x\in{\mathbb R}^N$ and $\sigma>0$. 
In the case of $1<p\le p_N$, inequality~\eqref{eq:1.2} is equivalent to 
$$
\sup_{x\in{\mathbb R}^N}\theta(B(x,\sqrt{T}))\le \gamma_1T^{\frac{N}{2}-\frac{1}{p-1}}. 
$$
In addition, if $p=p_N$, then
\begin{itemize}
  \item[(F2)] 
  there exists $\gamma_2=\gamma_2(N)>0$ such that
  $$
  \sup_{x\in{\mathbb R}^N}\theta(B(x,\sigma))\le\gamma_2\left[\log\left(e+\frac{\sqrt{T}}{\sigma}\right)\right]^{-\frac{N}{2}}
  \quad
  \mbox{for all $\sigma\in(0,\sqrt{T})$}.
  $$ 
\end{itemize}
See also \cites{ADi,BP} for properties~(F1) and (F2). 

On the other hand, 
let $u$ be a nonnegative solution to problem~\eqref{eq:P'} in ${\mathbb R}^N\times(0,T)$, that is, 
$u$ satisfies 
$$
u(x,t)=\int_{{\mathbb R}^N}\Gamma_N(x-y,t)\,d\eta(y)
+\int_0^t\int_{{\mathbb R}^N}\Gamma_N(x-y,t-s)u(y,s)^p\,dy\,ds
$$
for a.a.~$(x,t)\in{\mathbb R}^N\times(0,T)$. 
Then the initial data $\eta$ of the solution $u$ coincides with the initial trace of the solution~$u$ (see \cite{HI18}*{Theorem~1.2})
and it must satisfy properties~(F1) and~(F2) with~$\theta$ replaced by $\eta$.

Sufficient conditions on the initial data $\eta$ for the solvability of problem~\eqref{eq:P'} have been studied in many papers 
(see e.g. \cites{HI18, RS, BC, FHIL, FIoku, IKO, IKO2, KY, LRSV, LS, M, T, W2, W80}), 
and we have the following results.
\begin{itemize}
  \item[(F3)] Let $1<p<p_N$. Then there exists $\gamma_3=\gamma_3(N,p)>0$ such that, if
 \[
 \sup_{x\in{\bf R}^N}\eta(B(x,\sqrt{T}))\le \gamma_3T^{\frac{N}{2}-\frac{1}{p-1}}
 \]
 for some $T>0$, then problem~\eqref{eq:P'} possesses a solution in ${\bf R}^N\times(0,T)$. 
 (See e.g. \cite{HI18}*{Theorem~1.3}.)
\item[(F4)] Let $p\ge p_N$. Assume that 
\[
0\leq \eta(x)\leq 
\left\{
\begin{array}{ll}
 \gamma|x|^{-N}\displaystyle{\biggr[\log\left(e+\frac{1}{|x|}\right)\biggr]^{-\frac{N}{2}-1}} & \quad\mbox{if}\quad p=p_N,\vspace{7pt}\\
 \gamma|x|^{-\frac{2}{p-1}} & \quad\mbox{if}\quad p>p_N,\vspace{3pt}\\
\end{array}
\right.
\]
for some $\gamma>0$. 
Then there exists $\gamma_4=\gamma_4(N,p)>0$ such that
problem~\eqref{eq:P'} possesses a local-in-time solution if $\gamma\le \gamma_4$. 
(See \cite{HI18}*{Corollary~1.2}.)
\end{itemize}
The results in (F3) and (F4) show the optimality of necessary conditions given in~(F1) and (F2). 
Furthermore, 
the results in (F1), (F2), and (F4) imply
that the strength of the singularity  at the origin of the functions
\[
\eta(x)=
\left\{
\begin{array}{ll}
 |x|^{-\frac{2}{p-1}} & \quad\mbox{if}\quad p>p_N,\vspace{3pt}\\
 |x|^{-N}|\log|x||^{-\frac{N}{2}-1} & \quad\mbox{if}\quad p=p_N,
\end{array}
\right.
\]
is the critical threshold 
for the solvability of problem~\eqref{eq:P'}. 
We term such a singularity in the initial data an {\em optimal singularity} of initial data
for the solvability of problem~\eqref{eq:P'}. 
Recently, optimal singularities of initial data for the solvability 
were studied extensively for various nonlinear parabolic problems 
(see e.g. \cites{FI01,FI02,HI19,HT,HIT,HS,IKS,IS01,IS02}). 
However, these are not available to our problem~\eqref{eq:P} due to the zero Dirichlet boundary condition. 
In this paper we achieve our objectives~(O1)--(O4), and then we find optimal singularities 
of the initial data for the solvability of problem~(P), 
which are more complicated than those of problem~\eqref{eq:P'} 
and have the following three threshold cases:  
\[
{\rm (i)}\quad p=p_{N+1};\qquad {\rm (ii)}\quad p=p_N;\qquad {\rm (iii)}\quad p=2.
\]
See Section~5.3.
\subsection{Notation and definition of solutions}
We introduce some notation. 
As already said, unless otherwise stated, we set 
$\Omega:={\mathbb R}^{N-1}\times(0,\infty)$ if $N\ge 2$ and $\Omega:=(0,\infty)$ if $N=1$. 
We often identify $\partial\Omega={\mathbb R}^{N-1}$ if $N\ge 2$ and $\partial\Omega=\{0\}$ if $N=1$. 
For any $x=(x',x_N)\in\overline{\Omega}$, $r>0$, and $L>0$, let 
\begin{equation*}
\begin{array}{ll}
B_\Omega(x,r):=B(x,r)\cap\overline{\Omega},  & B'(x',r) := \{y'\in{{\mathbb R}}^{N-1} : |x'-y'| <r\}\subset\partial\Omega,\vspace{3pt}\\
\Omega_L:=\{(x',x_N)\in \Omega\,: x_N\ge L\}. &
\end{array}
\end{equation*}

We denote by ${\mathcal M}$ the set of nonnegative Radon measures on $\overline{\Omega}$, that is, 
$\mbox{{\rm supp}}\,\mu\subset\overline{\Omega}$ for $\mu\in {\mathcal M}$. 
For any $L^1_{\rm loc}(\overline{\Omega})$-function $\mu$,  
we often identify $d\mu=\mu (x)dx$ in ${\mathcal M}$.
For any $T\in(0,\infty]$, 
we set $Q_T:=\Omega\times(0,T)$. 
We denote by~${\mathcal L}(Q_T)$ the set of nonnegative, measurable, and finite almost everywhere functions in $Q_T$. 

For any $d=1,2,\dots$, let $\Gamma_d$ be the heat kernel in ${\mathbb R}^d\times(0,\infty)$, that is, 
\begin{equation}
\label{eq:1.3}
\Gamma_d(x,t):= (4\pi t)^{-\frac{d}{2}}\exp\left(-\frac{|x|^2}{4t}\right)
\quad\mbox{for}\quad (x,t)\in{\mathbb R}^d\times(0,\infty).
\end{equation}
Let $G$ be the Dirichlet heat kernel in $Q_\infty$, that is, 
\begin{equation}
\label{eq:1.4}
\begin{split}
G(x,y,t) := &\, \Gamma_{N-1}(x'-y',t)
 \left[\Gamma_1(x_N-y_N,t)- \Gamma_1(x_N+y_N,t)\right]\\
 = &\, \Gamma_N(x-y,t)
 \left(1-\exp\left(-\frac{x_Ny_N}{t}\right)\right)
\end{split}
\end{equation}
for $x=(x',x_N)$, $y=(y',y_N)\in\overline{\Omega}$, and $t>0$.
Then 
\begin{equation*}
\left\{
\begin{array}{ll}
 G(x,y,t)=G(y,x,t)\quad & \mbox{if}\quad (x,y,t)\in\overline{\Omega}\times\overline{\Omega}\times(0,\infty),\vspace{3pt}\\
 G(x,y,t)>0\quad & \mbox{if}\quad (x,y,t)\in\Omega\times\Omega\times(0,\infty),\vspace{3pt}\\
 G(x,y,t)=0\quad & \mbox{if}\quad (x,y,t)\in\partial\Omega\times\overline{\Omega}\times(0,\infty).
\end{array}
\right.
\end{equation*}
Furthermore, 
\begin{equation}
\label{eq:1.5}
\begin{split}
 & \lim_{y_N\to +0}y_N^{-1}G(x,y',y_N,t)=(\partial_{y_N}G)(x,y',0,t)\\
 & \qquad\quad
 =-2\Gamma_{N-1}(x'-y',t)\,\partial_{x_N}\Gamma_1(x_N,t)
=\frac{x_N}{t} \Gamma_N(x'-y',x_N,t)
\end{split}
\end{equation}
for $x=(x',x_N)\in\overline{\Omega}$, $y'\in{\mathbb R}^{N-1}$, and $t>0$. 
Define 
\begin{equation}
\label{eq:1.6}
K(x,y,t):=
\left\{
\begin{array}{ll}
y_N^{-1}G(x,y,t)\quad & \mbox{if}\quad y_N>0,\vspace{5pt}\\
(\partial_{y_N}G)(x,y,t)\quad & \mbox{if}\quad y_N=0,
\end{array}
\right.
\end{equation}
for $x\in\overline{\Omega}$, $y=(y',y_N)\in\overline{\Omega}$, and $t>0$. 
Then $K\in C(\overline{\Omega}\times\overline{\Omega}\times(0,\infty))$ and 
\begin{equation*}
\left\{
\begin{array}{ll}
 K(x,y,t)>0\quad & \mbox{if}\quad (x,y,t)\in\Omega\times\overline{\Omega}\times(0,\infty),\vspace{3pt}\\
 K(x,y,t)=0\quad & \mbox{if}\quad (x,y,t)\in\partial\Omega\times\overline{\Omega}\times(0,\infty).
\end{array}
\right.
\end{equation*}

We formulate definitions of solutions to problems~(E) and \eqref{eq:P}. 
\begin{definition}
\label{Definition:1.1}
Let $p>1$ and $T\in(0,\infty]$.
\begin{itemize}
  \item[{\rm (i)}] 
  We say that a function $u\in {\mathcal L}(Q_T)$ solves problem~{\rm (E)} in $Q_T$
  if $u$ satisfies
  \begin{equation}
  \label{eq:1.7}
  u(x,t)=\int_\Omega G(x,y,t-\tau)u(y,\tau)\,dy+\int_\tau^t\int_\Omega G(x,y,t-s)u(y,s)^p\,dy\,ds
  \end{equation}
  for a.a.~$(x,t)\in Q_T$ and a.a.~$\tau\in(0,T)$ with $\tau<t$.
  \item[{\rm (ii)}] 
  Let $\mu\in{\mathcal M}$. 
  We say that a function $u\in {\mathcal L}(Q_T)$ 
  is a solution to problem~\eqref{eq:P} in $Q_T$ if $u$ satisfies
  \begin{equation}
  \label{eq:1.8}
  u(x,t)=\int_{\overline{\Omega}}K(x,y,t)\,d\mu(y)+\int_0^t\int_\Omega G(x,y,t-s)u(y,s)^p\,dy\,ds
  \end{equation}
  for a.a.~$(x,t)\in Q_T$. 
  If $u\in {\mathcal L}(Q_T)$ satisfies \eqref{eq:1.8} with $``="$ replaced by $``\ge"$,
  then we say that $u$ is a supersolution to problem~\eqref{eq:P} in $Q_T$.
  \end{itemize}
\end{definition}
\begin{remark}
\label{Remark:1.1}
Let $\mu\in{\mathcal M}$. Set $\eta:=x_N^{-1}\mu$, which is a nonnegative Borel regular measure in $\Omega$. 
Let $u\in {\mathcal L}(Q_T)$ be a solution to problem~\eqref{eq:P} in $Q_T$, where $T>0$. 
It follows from Definition~{\rm\ref{Definition:1.1}~(ii)} that 
\begin{align*}
u(x,t)=\int_{\Omega}G(x,y,t)\,d\eta(y)+\int_{\partial\Omega}(\partial_{y_N}G)(x,y,t)\,d\mu\big|_{\partial\Omega}(y')+\int_0^t\int_\Omega G(x,y,t-s)u(y,s)^p\,dy\,ds
\end{align*}
for a.a.~$(x,t)\in Q_T$. 
Then the solution~$u$ can be regarded as a solution to problem~{\rm (E)} with the initial data
$$
u(x,0)=\eta-\mu\big|_{\partial\Omega}\otimes\delta_1'\quad\mbox{on}\quad\overline{\Omega},
$$ 
where $\delta_1'$ is the distributional derivative of the $1$-dimensional Dirac measure concentrated at the origin. 
Similarly,  the solution~$u$ can also be regarded as a solution to problem~{\rm (E)} with 
\begin{equation*}
\left\{
\begin{array}{ll}
u(x,0)=\eta\quad\mbox{in}\quad\Omega\qquad & \mbox{{\rm ({\it the initial data})}};\vspace{4pt}\\
u(x,t)=\mu\big|_{\partial\Omega}\otimes\delta_1\quad\mbox{on}\quad\partial\Omega\times[0,T)\quad & \mbox{{\rm ({\it the lateral boundary condition})}}.
\end{array}
\right.
\end{equation*}
\end{remark}
\subsection{Main results}
We state our main results on initial traces of solutions to problem~(E) and 
necessary conditions for the solvability of problem~\eqref{eq:P}. 
For our results on sufficient conditions and optimal singularities of the initial data, see Section~5. 
\vspace{3pt}

Theorem~\ref{Theorem:1.1} is concerned with the existence and the uniqueness of initial trace of $x_Nu(x,t)$ in ${\mathcal M}$. 
\begin{theorem}
\label{Theorem:1.1}
Let $N\geq1$, $p>1$, and $T>0$. 
\begin{itemize}
  \item[{\rm (i)}] 
  If $u\in {\mathcal L}(Q_T)$ solves problem~\eqref{eq:E},  
  then there exists a unique $\nu\in{\mathcal M}$ such that 
  \begin{equation}
  \label{eq:1.9}
  \underset{t\to +0}{\mbox{{\rm ess lim}}}
  \int_{\Omega} y_Nu(y,t)\phi(y)\,dy=\int_{\overline{\Omega}} \phi(y)\,d\nu(y)
  \end{equation}
  for all $\phi\in C_0({\mathbb R}^N)$. 
  Furthermore, $u$ is a solution to problem~\eqref{eq:P} with $\mu=\nu$ in $Q_T$.
  \item[{\rm (ii)}] 
  Let $\mu\in{\mathcal M}$. 
  If $u\in {\mathcal L}(Q_T)$ is a solution to problem~\eqref{eq:P} in $Q_T$, 
  then $u$ solves problem~{\rm (E)} in $Q_T$ and it satisfies \eqref{eq:1.9} with $\nu=\mu$.
\end{itemize}
\end{theorem}
In Theorem~\ref{Theorem:1.2} 
we obtain necessary conditions for the existence of local-in-time 
supersolutions to problem~\eqref{eq:P}. 
\begin{theorem}
\label{Theorem:1.2}
Let $N\geq1$, $p>1$, $T>0$, and $\mu\in{\mathcal M}$. 
Let $u\in {\mathcal L}(Q_T)$ be a supersolution to problem~\eqref{eq:P} in $Q_T$. 
There exists $\gamma_1=\gamma_1(N,p)>0$ such that
\begin{equation}
\label{eq:1.10}
\mu(B_\Omega(z,\sigma))\leq \gamma_1 
\sigma^{-\frac{2}{p-1}}\int_{B_\Omega(z,\sigma)}y_N\,dy
\end{equation}
for all $z\in\overline{\Omega}$ and $\sigma\in(0,\sqrt{T})$. 
  In addition,
  \begin{itemize}
  \item[{\rm (i)}]
  if $p=p_N$, then
  there exists $\gamma_2=\gamma_2(N)>0$ such that 
  \begin{equation}
  \label{eq:1.11}
  z_N^{-1}\mu(B_\Omega(z,\sigma))\le\gamma_2 \left[\log\left(e+\frac{\sqrt{T}}{\sigma}\right)\right]^{-\frac{N}{2}}
  \end{equation}
  for all $z=(z',z_N)\in\Omega_{3\sigma}$ and $\sigma\in(0,\sqrt{T})$;
  \item[{\rm (ii)}]
  if $p=p_{N+1}$, then there exists $\gamma_3=\gamma_3(N)>0$ such that
  \begin{equation}
  \label{eq:1.12}
  \mu(B_\Omega(z,\sigma))\le\gamma_3
  \left[\log\left(e+\frac{\sqrt{T}}{\sigma}\right)\right]^{-\frac{N+1}{2}}
  \end{equation}
  for all $z\in\partial\Omega$ and $\sigma\in(0,\sqrt{T})$;
  \item[{\rm (iii)}]
  if $p\ge 2$, then $\mu(\partial\Omega)= 0$. 
  \end{itemize}
\end{theorem}
\begin{remark}
\label{Remark:1.2}
Assertion~{\rm (iii)} with $p>2$ follows from \eqref{eq:1.10}.
Indeed, by \eqref{eq:1.10}, 
applying covering theorems {\rm({\it see e.g. \cite{EG}*{Section~1.5}})}, 
we have
$$
\mu(B'(x',\sqrt{T})\times[0,\sigma))\le C\sigma^{-\frac{2}{p-1}}\int_{B'(x',\sqrt{T}+\sigma)\times[0,2\sigma)}y_N\,dy
\le CT^{\frac{N-1}{2}}\sigma^{2-\frac{2}{p-1}}
$$
for all $x'\in\partial\Omega$ and $\sigma\in(0,\sqrt{T})$. 
Then $\mu(B'(x',\sqrt{T})\times\{0\})=0$ for all $x'\in\partial\Omega$, that is, $\mu(\partial\Omega)=0$. 
\end{remark}
As a corollary of Theorems~\ref{Theorem:1.1} and \ref{Theorem:1.2}, 
we have:
\begin{corollary}
\label{Corollary:1.1} 
Let $p\ge 2$ and $T>0$. 
Let $u\in {\mathcal L}(Q_T)$ solve problem~{\rm (E)} in $Q_T$. 
Then there exists a unique $\nu\in{\mathcal M}$ such that 
$u$ is a solution to problem~{\rm (P)} with $\mu=\nu$ in $Q_T$ and $\nu(\partial\Omega)=0$, that is, $u$ satisfies 
$$
u(x,t)=\int_\Omega G(x,y,t)y_N^{-1}\,d\nu(y)+\int_0^t\int_\Omega G(x,y,t-s)u(y,s)^p\,dy\,ds
\quad \mbox{for a.a.~$(x,t)\in Q_T$}.
$$
\end{corollary}
Furthermore, as an application of Theorem~\ref{Theorem:1.2}, 
we obtain the following result on the blow-up rate of solutions to problem~\eqref{eq:P}.
\begin{corollary}
\label{Corollary:1.2} 
Let $p>1$ and $T>0$. 
Let $u\in {\mathcal L}(Q_T)$ solve problem~{\rm (E)} in $Q_T$. 
\begin{itemize}
  \item[{\rm (i)}] 
  There exists $C_1=C_1(N,p)>0$ such that 
  \begin{equation*}
  \int_{B(x,\sqrt{T-t})}u(y,t)\,dy\le C_1(T-t)^{\frac{N}{2}-\frac{1}{p-1}}
  \end{equation*}
  for all $x=(x',x_N)\in\Omega$ with $x_N\ge\sqrt{T-t}$ and a.a.~$t\in(0,T)$. 
  \item[{\rm (ii)}]
  There exists $C_2=C_2(N,p)>0$ such that 
  \begin{equation*}
  \int_{B_\Omega(x,\sqrt{T-t})}y_Nu(y,t)\,dy\le C_2(T-t)^{\frac{N+1}{2}-\frac{1}{p-1}}
  \end{equation*}
  for all $x\in\partial\Omega$ and a.a.~$t\in(0,T)$. 
\end{itemize}
\end{corollary}
See e.g. \cite{QS}*{Section~23} for further information on blow-up rates of solutions. 
\vspace{3pt}

We explain the outline of the proofs of Theorems~\ref{Theorem:1.1} and \ref{Theorem:1.2}. 
For the proof of Theorem~\ref{Theorem:1.2}, 
we prepare some lower estimates of the Dirichlet heat kernel $G$ (see Lemmas~\ref{Lemma:3.1} and \ref{Lemma:3.2}). 
Then, applying the arguments in the proofs of \cite{HI18}*{Theorem~1.1} and \cite{FHIL}*{Theorem~3.1}, 
we obtain inequalities~\eqref{eq:1.10} and \eqref{eq:1.11} for balls having a positive distance from the boundary~$\partial\Omega$ 
(see Proposition~\ref{Proposition:3.1}). 
Furthermore, combining lower estimates of $G$ and a covering lemma, 
we prove inequalities \eqref{eq:1.10} and \eqref{eq:1.11}. 
In addition, we prove inequality~\eqref{eq:1.12} by using some decay estimates 
of the integral kernels~$G$ and $K$ near the boundary~$\partial\Omega$. 
Next, 
we show that the existence of solution $u$ to problem~\eqref{eq:P} in $Q_T$ implies that
the function
\[
U(x_N,t):= \int_{{\mathbb R}^{N-1}} \Gamma_{N-1}(x',T)u(x+z,t) \,dx'
\]
is a supersolution to problem~\eqref{eq:P} with $N=1$. 
Then, by Theorem~\ref{Theorem:1.2}~(ii) and \eqref{eq:1.10} we prove that $\mu(\partial\Omega)= 0$ if $p\ge 2$, 
and complete the proof of Theorem~\ref{Theorem:1.2}.

Theorem~\ref{Theorem:1.1} is proved by Theorem~\ref{Theorem:1.2}. 
For any function~$u\in {\mathcal L}(Q_T)$ solving problem~(E) in $Q_T$, 
by Theorem~\ref{Theorem:1.2} we obtain uniform local estimates of $x_Nu(x,t)$, 
and prove the existence and the uniqueness of the initial trace of $x_Nu(x,t)$.
Then we modify the arguments in the proof of \cite{HI18}*{Theorem~1.2} to prove Theorem~\ref{Theorem:1.1}. 
Inequality~\eqref{eq:1.10} and the Besicovitch covering lemma are used effectively in the proof of Theorem~\ref{Theorem:1.1}. 
\vspace{3pt}

The rest of this paper is organized as follows. 
In Section~2 we prepare some preliminary lemmas. 
In Sections~3 and 4 we prove Theorems~\ref{Theorem:1.2} and \ref{Theorem:1.1}, respectively. 
The proofs of Corollaries~\ref{Corollary:1.1} and \ref{Corollary:1.2} are also given in Section~4. 
In Section~5 we modify the arguments in \cites{HI18, RS} to 
obtain sufficient conditions for the solvability of problem~\eqref{eq:P}. 
Furthermore, 
we combine our necessary conditions and sufficient conditions 
to find optimal singularities of the initial data for the solvability of problem~\eqref{eq:P}.
\section{Preliminaries}
In what follows we will use $C$ to denote generic positive constants. 
The letter $C$ may take different values within a calculation.
We first prove the following covering lemma.  
\begin{lemma}
\label{Lemma:2.1}
Let $N\ge 1$ and $\delta\in(0,1)$. 
Then there exists $m\in\{1,2,\dots\}$ with the following properties. 
\begin{itemize}
  \item[{\rm (i)}] 
  For any $z\in{\mathbb R}^N$ and $r>0$, 
  there exists $\{z_i\}_{i=1}^m\subset{\mathbb R}^N$ such that 
  $$
  B(z,r)\subset\bigcup_{i=1}^m B(z_i,\delta r).
  $$ 
  \item[{\rm (ii)}] 
  For any $z\in{\mathbb R}^N$ and $r>0$, 
  there exists $\{\overline{z}_i\}_{i=1}^m\subset B_\Omega(z,2r)$ such that 
  $$
  B_\Omega(z,r)\subset\bigcup_{i=1}^m B_\Omega(\overline{z}_i,\delta r).
  $$ 
\end{itemize}
\end{lemma}
{\bf Proof.}
Let $\delta\in(0,1)$, $r>0$, and $z\in{\mathbb R}^N$. 
We find $m\in\{1,2,\dots\}$ and $\{z_i\}_{i=1}^m\subset B(0,1)$ such that 
$B(0,1)\subset\cup_{i=1}^m B(z_i,\delta)$, so that
\begin{equation}
\label{eq:2.1}
B(z,r)\subset\bigcup_{i=1}^m B(z+rz_i,\delta r),
\end{equation}
which implies assertion~(i). 
Similarly, 
by \eqref{eq:2.1} we find $m'\in\{1,2,\dots\}$ and $\{\tilde{z}_i\}_{i=1}^{m'}\subset B_\Omega(0,1)$ such that 
$B_\Omega(0,1)\subset\cup_{i=1}^{m'} B_\Omega(\tilde{z}_i,\delta/2)$, so that
\begin{equation}
\label{eq:2.2}
B_\Omega(z,r)\subset\bigcup_{i=1}^{m'} B_\Omega(z+r \tilde{z}_i,\delta r/2).
\end{equation}
Set $\overline{z}_i:=z+r\tilde{z}_i$ if $z+r\tilde{z}_i\in\overline{\Omega}$ 
and $\overline{z}_i:=(z'+r \tilde{z}_i',0)$ if $z+r\tilde{z}_i\not\in\overline{\Omega}$. 
Then
\[
\overline{z}_i\in B_\Omega(z,2r),\qquad
B_\Omega(z+r \tilde{z}_i,\delta r/2)\subset B_\Omega(\overline{z}_i,\delta r)
\quad\mbox{if}\quad B_\Omega(z+r\tilde{z}_i,\delta r/2)\not=\emptyset.
\]
This together with \eqref{eq:2.2} implies that 
\[
B_\Omega(z,r)\subset\bigcup_{i=1}^{m'} B_\Omega(\overline{z}_i,\delta r).
\]
Then assertion~(ii) follows, and the proof is complete.
$\Box$
\vspace{5pt}

Next, we state two lemmas on the integral kernels $\Gamma_d$, $G$, and $K$. 
\begin{lemma}
\label{Lemma:2.2}
{\rm (i)} 
There exists $C_1>0$ such that 
$$
\int_{{\mathbb R}^N}\Gamma_N(x-y,t)\,d\mu(y)\le C_1t^{-\frac{N}{2}}\sup_{z\in{\mathbb R}^N}\mu(B(z,\sqrt{t}))
$$
for all nonnegative Radon measure $\mu$ on ${\mathbb R}^N$ and $(x,t)\in{\mathbb R}^N\times(0,\infty)$.
\vspace{3pt}
\newline
{\rm (ii)} There exists $C_2>0$ such that 
\begin{equation}
\label{eq:2.3}
K(x,y,t)\le C_2\frac{x_N}{(x_N+\sqrt{t})(y_N+\sqrt{t})}\Gamma_N(x-y,2t)
\end{equation}
for all $(x,y,t)\in\Omega\times\overline{\Omega}\times(0,\infty)$. 
Furthermore, there exists $C_3>0$ such that 
\begin{equation}
\label{eq:2.4}
\int_{\overline{\Omega}}K(x,y,t)\,d\mu(y)\le C_3t^{-\frac{N}{2}}\sup_{z\in{\overline{\Omega}}}
\int_{B_\Omega(z,\sqrt{t})}\frac{d\mu(y)}{y_N+\sqrt{t}}
\end{equation}
for all $\mu\in{\mathcal M}$ and $(x,t)\in Q_\infty$. 
\end{lemma}
{\bf Proof.} 
Assertion~(i) follows from \cite{HI18}*{Lemma~2.1}. 
We prove assertion~(ii). 
It follows that 
\[
1-e^{-ab} \leq \min\{1,ab\} \leq \min\{1,a\} \min\{1,b\} (1+|a-b|) 
\leq \frac{4 ab(1+|a-b|)}{(1+a)(1+b)} 
\]
for all $a$, $b>0$ (see e.g. \cite{MMZ09}*{Section~1.1}). 
Then, by \eqref{eq:1.4} we have
\begin{equation*}
\begin{split}
G(x,y,t) & =\Gamma_N(x-y,t)
 \left(1-\exp\left(-\frac{x_Ny_N}{t}\right)\right)\\
  & \le \Gamma_N(x-y,t)\frac{4x_N y_N}{(x_N+\sqrt{t})(y_N+\sqrt{t})}
  \left( 1+ \frac{|x_N-y_N|}{\sqrt{t}} \right)\\
  & = (4\pi t)^{-\frac{N}{2}} e^{-\frac{|x'-y'|^2}{4t}} 
  e^{-\frac{|x_N-y_N|^2}{8t}}
  \frac{y_N}{y_N+\sqrt{t}}
  \left( 1+ \frac{|x_N-y_N|}{\sqrt{t}} \right) 
  e^{-\frac{|x_N-y_N|^2}{8t}} 
  \frac{4x_N }{x_N+\sqrt{t}} \\
  & \le C\frac{x_N}{x_N+\sqrt{t}}\frac{y_N}{y_N+\sqrt{t}}\Gamma_N(x-y,2t)
\end{split}
\end{equation*}
for all $x=(x',x_N)$, $y=(y',y_N)\in\Omega$, and $t>0$. 
This together with \eqref{eq:1.5} and \eqref{eq:1.6} implies~\eqref{eq:2.3}. 
Furthermore, we have
$$
\int_{\overline{\Omega}}K(x,y,t)\,d\mu(y)\le C\int_{\overline{\Omega}}\Gamma_N(x-y,2t)\frac{d\mu(y)}{y_N+\sqrt{t}}
\quad\mbox{for $(x,t)\in Q_\infty$},
$$
which together with Lemma~\ref{Lemma:2.1} and assertion~(i) implies \eqref{eq:2.4}. 
Thus Lemma~\ref{Lemma:2.2} follows. 
$\Box$

\begin{lemma}
\label{Lemma:2.3}
The integral kernels $K$ and $G$ satisfy 
\begin{eqnarray}
\label{eq:2.5}
 & & \int_\Omega K(x,y,t)\,dx=(\pi t)^{-\frac{1}{2}}\quad\mbox{for}\,\,\,(y,t)\in\partial\Omega\times(0,\infty),\\
\label{eq:2.6}
 & & \int_\Omega G(z,x,s)K(x,y,t)  dx= K(z,y,t+s)\quad\mbox{for}\,\,\,(z,y,t,s)\in\Omega\times\overline{\Omega}\times(0,\infty)^2. 
\end{eqnarray}
\end{lemma}
{\bf Proof.}
Let $y=(y',0)\in\partial\Omega$. 
It follows from~\eqref{eq:1.5} that 
$$
\int_\Omega K(x,y,t)\,dx =\int_{{\mathbb R}^{N-1}}\Gamma_{N-1} (x'-y',t)\,dx'\int_0^\infty (4\pi t)^{-\frac{1}{2}}\frac{x_N}{t}\exp\left(-\frac{x_N^2}{4t}\right)\,dx_N
 =(\pi t)^{-\frac{1}{2}}, 
$$
which implies~\eqref{eq:2.5}. 
On the other hand, since 
\begin{equation*}
\begin{split}
&\int_0^\infty [\Gamma_1(x_N-y_N,t)-\Gamma_1(x_N+y_N,t)] 
(\partial_{y_N}\Gamma_1)(y_N,s) \, dy_N \\
&=-\int_0^\infty 
\partial_{ y_N} [\Gamma_1(x_N-y_N,t)-\Gamma_1(x_N+y_N,t)] \cdot \Gamma_1(y_N,s) \, dy_N\\
&=\partial_{ x_N}\int_0^\infty 
[\Gamma_1(x_N-y_N,t)+\Gamma_1(x_N+y_N,t)] \Gamma_1(y_N,s) \, dy_N\\
&=\partial_{ x_N}\int_{-\infty}^\infty \Gamma_1(x_N-y_N,t)\Gamma_1(y_N,s) \, dy_N
=(\partial_{x_N} \Gamma_1) (x_N,t+s) 
\end{split}
\end{equation*}
for $x_N\in(0,\infty)$ and $t$, $s\in(0,\infty)$, by~\eqref{eq:1.4} and \eqref{eq:1.5} 
we have 
\begin{equation*}
\begin{split}
 & \int_\Omega G(z,x,s) K(x,y,t)  dx\\
 & =-2\int_{{\mathbb R}^{N-1}}\Gamma_{N-1}(z'-x',s) \Gamma_{N-1}(x'-y',t)  \,dx'\\
 & \qquad\quad
  \times\int_0^\infty 
 [\Gamma_1(z_N-x_N,s) -\Gamma_1(z_N+x_N,s)] 
  (\partial_{x_N}\Gamma_{1})(x_N,t) \,dx_N \\
  & =-2\Gamma_{N-1}(z'-y',t+s) 
  (\partial_{x_N} \Gamma_1)(z_N,t+s)\\
  & =\frac{z_N}{t+s}  \Gamma_{N-1}(z'-y',t+s)\Gamma_{1}(z_N,t+s)
  =\frac{z_N}{t+s}  \Gamma_N(z-y,t+s)
  =K(z,y,t+s)
\end{split}
\end{equation*}
for $z\in \Omega$ and $t, s\in(0,\infty)$.
Thus relation~\eqref{eq:2.6} holds for $y\in\partial\Omega$. 
In the case of $y=(y',y_N)\in\Omega$, we have
$$
\int_\Omega G(z,x,s)K(x,y,t)  dx=y_N^{-1}\int_\Omega G(z,x,s)G(x,y,t)  dx
=y_N^{-1}G(z,y,t+s)=K(z,y,t+s)
$$
for $z\in \Omega$ and $t, s\in(0,\infty)$. 
Thus relation~\eqref{eq:2.6} holds for $y\in\Omega$, 
and the proof is complete.
$\Box$\vspace{5pt}

Next, we state preliminary lemmas on solutions to problem~\eqref{eq:P}. 
Lemmas~\ref{Lemma:2.4} and \ref{Lemma:2.5} follow from the nonnegativity of the integral kernels $G$ and $K$. 

\begin{lemma}
\label{Lemma:2.4}
Let $u\in {\mathcal L}(Q_T)$ be a solution to problem~\eqref{eq:P} in $Q_T$, where $T>0$.  
Then, for a.a.~$\tau\in(0,T)$, the function 
$u_\tau(x,t):=u(x,t+\tau)$
is a solution to problem~\eqref{eq:P} 
with $\mu=x_Nu(x,\tau)$ in $Q_{T-\tau}$. 
Furthermore, for a.a.~$\tau\in(0,T)$, 
$$
\int_\Omega G(x,y,t)u(y,\tau)\,dy<\infty
\quad\mbox{for a.a.~$(x,t)\in Q_{T-\tau}$}.
$$
\end{lemma}
{\bf Proof.}
It follows from Definition~\ref{Definition:1.1}, \eqref{eq:2.6}, 
and Fubini's theorem that 
\begin{equation*}
\begin{split}
 & \int_\Omega G(x,y,t)u(y,\tau)\,dy+\int_0^t\int_\Omega G(x,y,t-s)u(y,s+\tau)^p\,dy\,ds\\
 & =\int_\Omega\int_{\overline{\Omega}}G(x,y,t)K(y,z,\tau)\,d\mu(z)\,dy
 +\int_0^\tau\int_\Omega\int_\Omega G(x,y,t)G(y,z,\tau-s)u(z,s)^p\,dz\,dy\,ds\\
 & \qquad\quad+\int_\tau^{t+\tau}\int_\Omega G(x,y,t+\tau-s)u(y,s)^p\,dy\,ds\\
 & =\int_{\overline{\Omega}}K(x,z,t+\tau)\,d\mu(z)+\int_0^{t+\tau}\int_\Omega G(x,z,t+\tau-s)u(z,s)^p\,dz\,ds
 =u(x,t+\tau)<\infty
\end{split}
\end{equation*}
for a.a.~$(x,t)\in Q_{T-\tau}$ and a.a.~$\tau\in(0,T)$. 
Then Lemma~\ref{Lemma:2.4} follows.
$\Box$
\begin{lemma}
\label{Lemma:2.5}
Assume that there exists a supersolution~$v$ to problem~\eqref{eq:P} in $Q_T$. 
Then problem~{\rm (P)} possesses a solution~$u$ to problem~\eqref{eq:P} in $Q_T$ such that $u\le v$ in $Q_T$.
\end{lemma}
{\bf Proof.}
By the same argument as in the proof of \cite{HI18}*{Lemma~2.2} we find a solution~$u\in {\mathcal L}(Q_T)$ 
satisfying $u(x,t)\le v(x,t)$ for a.a.~$(x,t)\in Q_T$. 
Thus Lemma~\ref{Lemma:2.5} follows.
$\Box$\vspace{5pt}

At the end of this section we prepare a lemma on an integral inequality. 
The idea of using this kind of lemma is due to \cite{LS02}. (See also the proof of \cite{FHIL}*{Theorem~3.1}.)
\begin{lemma}
\label{Lemma:2.6}
Let $\zeta$ be a nonnegative measurable function in $(0,T)$, where $T>0$. 
Assume that 
\begin{equation}
\label{eq:2.7}
\infty>\zeta(t)\ge c_1+c_2\int_{t_*}^t s^{-\alpha}\zeta(s)^\beta\,ds
\quad \mbox{for a.a.~$t\in(t_*,T)$},
\end{equation}
where $c_1$, $c_2>0$, $\alpha\ge 0$, $\beta>1$, and $t_*\in(0,T/2)$. 
Then there exists $C=C(\alpha,\beta)>0$ such that 
\begin{equation}
\label{eq:2.8}
c_1\le C c_2^{-\frac{1}{\beta-1}}t_*^{\frac{\alpha-1}{\beta-1}}.
\end{equation}
In addition, 
if $\alpha=1$, then 
\begin{equation}
\label{eq:2.9}
c_1\le (c_2(\beta-1))^{-\frac{1}{\beta-1}}\left[\log\frac{T}{2t_*}\right]^{-\frac{1}{\beta-1}}.
\end{equation}
\end{lemma}
{\bf Proof.}
Let $\eta$ be a solution to the initial value problem 
$\eta'(t)=c_2t^{-\alpha}\eta(t)^\beta$ with $\eta(t_*)=c_1$. 
We observe from \eqref{eq:2.7} 
that $\eta$ exists in $[t_*,T)$.
Since $t_*<T/2$, 
we see that 
$$
\int_{\eta(t_*)}^\infty \xi^{-\beta}\,d\xi
\ge \int_{\eta(t_*)}^{\eta(2t_*)} \xi^{-\beta}\,d\xi=c_2\int_{t_*}^{2t_*} s^{-\alpha}\,ds
\ge Cc_2t_*^{-\alpha+1},
$$
so that 
\[
\frac{1}{\beta-1}c_1^{-\beta+1}\ge Cc_2t_*^{-\alpha+1}.
\]
This implies \eqref{eq:2.8}. Furthermore, if $\alpha=1$, then 
\[
\frac{1}{\beta-1}c_1^{-\beta+1}\ge\int_{\eta(t_*)}^\infty \xi^{-\beta}\,d\xi
\ge \int_{\eta(t_*)}^{\eta(T/2)} \xi^{-\beta}\,d\xi
=c_2\int_{t_*}^{T/2} s^{-1}\,ds
\ge c_2\log\frac{T}{2t_*},
\]
so that \eqref{eq:2.9} holds. 
Thus Lemma~\ref{Lemma:2.6} follows.
$\Box$
\section{Proof of Theorem~\ref{Theorem:1.2}}
In this section we study necessary conditions for the solvability of problem~\eqref{eq:P}, 
and prove Theorem~\ref{Theorem:1.2}. 
We first modify the arguments in \cite{HI18} to prove the following proposition. 
\begin{proposition}
\label{Proposition:3.1} 
Assume that there exists a supersolution to problem~\eqref{eq:P}
in $Q_T$, where $T>0$. 
Then there exists $C_1>0$ such that 
\begin{equation}
\label{eq:3.1}
z_N^{-1}\mu(B(z,\sigma))\le C_1\sigma^{N-\frac{2}{p-1}}
\end{equation}
for all $z=(z',z_N)\in\Omega_{\sqrt{T}}$ and $\sigma\in(0,\sqrt{T}/16)$.
Furthermore, if $p=p_N$, there exists $C_2>0$ such that 
\begin{equation}
\label{eq:3.2}
z_N^{-1}\mu(B(z,\sigma))\le 
C_2\left[\log\left(e+\frac{\sqrt{T}}{\sigma}\right)\right]^{-\frac{N}{2}}
\end{equation}
for all $z=(z',z_N)\in\Omega_{\sqrt{T}}$ and $\sigma\in(0,\sqrt{T}/16)$.
\end{proposition}
In order to prove Proposition~\ref{Proposition:3.1}, 
we prepare the following two lemmas on integral kernels.
\begin{lemma}
\label{Lemma:3.1}
For any $\rho\in(0,1)$, 
there exists $C>0$ such that
\[
\int_{B(z,\sigma)} K(z,y,\sigma^2) \,d\mu(y)\ge C\sigma^{-N}z_N^{-1}\mu(B(z,\sigma))
\]
for all $\mu\in{\mathcal M}$, $z\in\Omega_{\sqrt{T}}$, $\sigma\in(0,\rho \sqrt{T})$, and $T>0$. 
\end{lemma}
{\bf Proof.}
Let $\rho\in(0,1)$, $\sigma\in(0,\rho \sqrt{T})$, $z=(z',z_N)\in \Omega_{\sqrt{T}}$, and $y=(y',y_N)\in B(z,\sigma)$. 
Since $z_N y_N\ge \sqrt{T}(z_N-\sigma)\ge(1-\rho)T$, 
by \eqref{eq:1.4} and \eqref{eq:1.6} we have
\begin{equation*}
\begin{split}
K(z,y,\sigma^2) & =y_N^{-1}G(z,y,\sigma^2)=y_N^{-1}\Gamma_N(z-y,\sigma^2)
\left(1-\exp\left(-\frac{z_Ny_N}{\sigma^2}\right)\right)\\
 & \ge Cz_N^{-1}\Gamma_N(z-y,\sigma^2)\left(1-\exp\left(-\frac{1-\rho}{\rho^2}\right)\right)
\ge Cz_N^{-1}\Gamma_N(z-y,\sigma^2).
\end{split}
\end{equation*}
This implies the desired inequality. The proof is complete. 
$\Box$
\begin{lemma}
\label{Lemma:3.2}
{\rm (i)} 
Let $d=1,2,\dots$. Then 
\[
\Gamma_d(x,2t-s)\ge\left(\frac{s}{2t}\right)^{\frac{d}{2}}\Gamma_d(x,s)
\]
for all $x\in{\mathbb R}^d$ and $t$, $s>0$ with $s<t$.
\vspace{3pt}
\newline
{\rm (ii)} 
There exists $C>0$ such that 
\[
G(z,y,2t-s) \geq C\left(\frac{s}{2t} \right)^\frac{N}{2} G(z,y,s)
\]
for all $z\in\Omega_{\sqrt{T}}$, $y\in\Omega$,  
$s, t\in (0,T/32)$ with $s<t$, and $T>0$.
\end{lemma}
{\bf Proof.}
We first prove assertion~(i).
Let $0<s<t$. It follows that 
\begin{equation*}
\begin{split}
\Gamma_d(x,2t-s)
 & =(4\pi s)^{-\frac{d}{2}} \left(\frac{s}{2t-s}\right)^{\frac{d}{2}}\exp\left(-\frac{|x|^2}{4(2t-s)}\right)\\
 & \ge(4\pi s)^{-\frac{d}{2}} \left(\frac{s}{2t}\right)^{\frac{d}{2}}\exp\left(-\frac{|x|^2}{4s}\right)
 =\left(\frac{s}{2t}\right)^{\frac{d}{2}} \Gamma_d(x,s)
\end{split}
\end{equation*}
for all $x\in{\mathbb R}^d$. This implies assertion~(i). 

Next, we prove assertion~(ii). 
Let $z\in\Omega_{\sqrt{T}}$, $y\in\Omega$, $s, t\in (0,T/32)$ 
with $s<t$, and $\sigma>0$.
If $y_N\geq z_N/2$, then, 
by \eqref{eq:1.4} and assertion~(i) we have 
\begin{equation*}
\begin{split}
G(z,y,2t-s) & =\Gamma_N(z-y,2t-s)
\left(1-\exp\left(-\frac{z_Ny_N}{2t-s}\right)\right)\\
 & \ge \left(\frac{s}{2t}\right)^{\frac{N}{2}}\Gamma_N(z-y,s)\left(1-\exp\left(-\frac{z_N^2}{4t}\right)\right)\\
 & \ge C\left(\frac{s}{2t}\right)^{\frac{N}{2}}\Gamma_N(z-y,s)
 \ge C\left(\frac{s}{2t}\right)^{\frac{N}{2}}G(z,y,s).
\end{split}
\end{equation*}
This implies assertion~(ii) in the case of $y_N\geq z_N/2$. 

Consider the case of $y_N<z_N/2$. 
Set 
\[
f(\tau) :=(4\pi\tau)^{\frac{1}{2}}\left[
\Gamma_1(z_N-y_N,\tau)-\Gamma_1 (z_N+y_N,\tau)\right]
\quad\mbox{for $\tau\in(0,2t)$}.
\]
Then 
\begin{equation*}
\begin{split}
f'(\tau) &
 =\frac{|z_N-y_N|^2}{4\tau^2}\exp\left(-\frac{|z_N-y_N|^2}{4\tau}\right)-\frac{|z_N+y_N|^2}{4\tau^2}\exp\left(-\frac{|z_N+y_N|^2}{4\tau}\right)\\
 & =\frac{1}{4\tau^2}\exp\left(-\frac{|z_N+y_N|^2}{4\tau}\right)\\
 & \qquad\times
 \left[|z_N-y_N|^2\exp\left(-\frac{|z_N-y_N|^2}{4\tau}+\frac{|z_N+y_N|^2}{4\tau}\right)-|z_N+y_N|^2\right]\\
 & =\frac{1}{4\tau^2}\exp\left(-\frac{|z_N+y_N|^2}{4\tau}\right)
\left[\exp\left(\frac{z_Ny_N}{\tau}\right)|z_N-y_N|^2-|z_N+y_N|^2\right]\\
& \ge\frac{1}{4\tau^2}\exp\left(-\frac{|z_N+y_N|^2}{4\tau}\right)
\left[\left(1+\frac{z_Ny_N}{\tau}\right)|z_N-y_N|^2-|z_N+y_N|^2\right]\\
 & =\frac{1}{4\tau^2}\exp\left(-\frac{|z_N+y_N|^2}{4\tau}\right)
z_Ny_N\left(-4+\frac{|z_N-y_N|^2}{\tau}\right).
\end{split}
\end{equation*}
Since $0<y_N<z_N/2$, $z_N\ge \sqrt{T}$, and $\tau\le 2t<T/16$, 
we have 
\[
\frac{|z_N-y_N|^2}{\tau}>\left(\frac{z_N}{2}\right)^2\frac{16}{T}\ge 4.
\]
Then we see that $f'\ge 0$ in $(0,2t)$. 
Since $2t>2t-s>s$, 
we observe that
\begin{equation*}
\begin{split}
 & \Gamma_1(z_N-y_N,2t-s)-\Gamma_1(z_N+y_N,2t-s)\\
 &=(4\pi(2t-s))^{-\frac{1}{2}}f(2t-s)\ge (4\pi(2t-s))^{-\frac{1}{2}}f(s)\\
 & =\left(\frac{s}{2t-s}\right)^{\frac{1}{2}}[\Gamma_1(z_N-y_N,s)-\Gamma_1(z_N+y_N,s)]\\
 & \ge\left(\frac{s}{2t}\right)^{\frac{1}{2}}[\Gamma_1(z_N-y_N,s)-\Gamma_1(z_N+y_N,s)]. 
\end{split}
\end{equation*}
Therefore, by \eqref{eq:1.4} and assertion~(i) we obtain
\begin{equation*}
\begin{split}
 & G(z,y,2t-s)=\Gamma_{N-1}(z'-y',2t-s)[\Gamma_1(z_N-y_N,2t-s)-\Gamma_1(z_N+y_N,2t-s)]\\
 & \qquad
 \ge \left(\frac{s}{2t}\right)^{\frac{N}{2}}
\Gamma_{N-1}(z'-y',s)[\Gamma_1(z_N-y_N,s)-\Gamma_1(z_N+y_N,s)]
=\left(\frac{s}{2t}\right)^{\frac{N}{2}}G(z,y,s). 
\end{split}
\end{equation*}
This implies assertion~(ii) in the case of $y_N<z_N/2$. 
Thus assertion~(ii) follows, and the proof is complete.  
$\Box$
\vspace{3pt}

\noindent
{\bf Proof of Proposition~\ref{Proposition:3.1}.} 
Let $u$ be a supersolution to problem~\eqref{eq:P} in $Q_T$, where $T>0$. 
Let $\sigma\in(0,\sqrt{T}/16)$ and $z=(z',z_N)\in\Omega_{\sqrt{T}}$. 
It follows from Lemma~\ref{Lemma:2.3} and Lemma~\ref{Lemma:3.2}~(ii) that 
\begin{equation*}
\begin{split}
 \int_\Omega G(z,x,t)u(x,t)\,dx & \ge\int_{\overline{\Omega}} \int_\Omega G(z,x,t)K(x,y,t)\,dx \,d\mu(y)\\
 & \qquad\quad
 +\int_0^t\int_\Omega\int_\Omega 
 G(z,x,t)G(x,y,t-s)u(y,s)^p\,dx \,dy\,ds\\
 & \ge\int_{\overline{\Omega}} K(z,y,2t)\,d\mu(y)
 +\int_{\sigma^2}^t\int_\Omega G(z,y,2t-s)u(y,s)^p\,dy\,ds\\
 & \ge \int_{\overline{\Omega}} K(z,y,2t)\,d\mu(y)
 +Ct^{-\frac{N}{2}}\int_{\sigma^2}^t\int_\Omega s^{\frac{N}{2}}G(z,y,s)u(y,s)^p\,dy\,ds
\end{split}
\end{equation*}
for a.a.~$t\in(\sigma^2,T/32)$.
Furthermore, 
H\"older's inequality implies that 
\begin{equation*}
\begin{split}
	\int_\Omega G(z,y,s)u(y,s)\,dy
	 & \le\left(\int_\Omega G(z,y,s)\,dy\right)^{1-\frac{1}{p}}\left(\int_\Omega 	G(z,y,s)u(y,s)^p\,dy\right)^{\frac{1}{p}}\\
	 & \le\left(\int_\Omega G(z,y,s)u(y,s)^p\,dy\right)^{\frac{1}{p}}
\end{split}
\end{equation*}
for all $s>0$. 
Then we obtain 
\begin{equation}
\label{eq:3.3}
\begin{split}
 & \int_\Omega G(z,x,t)u(x,t)\,dx\\
 & \ge  \int_{\overline{\Omega}} K(z,y,2t)\,d\mu(y)
 +Ct^{-\frac{N}{2}}\int_{\sigma^2}^t s^\frac{N}{2}
\left(\int_\Omega G(z,y,s)u(y,s)\,dy\right)^p\,ds
\end{split}
\end{equation}
for a.a.~$t\in(\sigma^2,T/32)$. 
In addition, 
Lemma~\ref{Lemma:3.1} implies that 
\begin{equation}
\label{eq:3.4}
t^{\frac{N}{2}}\int_{\overline{\Omega}} K(z,y,2t)\,d\mu(y)
\ge Ct^{\frac{N}{2}}(2t)^{-\frac{N}{2}}z_N^{-1}\mu(B(z,\sqrt{2t}))
\ge Cz_N^{-1}\mu(B(z,\sigma))
\end{equation}
for all $t\in(\sigma^2,T/32)$. 
Therefore, setting 
\[
W(t):=t^{\frac{N}{2}}\int_\Omega G(z,y,t)u(y,t)\,dy, 
\]
by \eqref{eq:3.3} and \eqref{eq:3.4} we obtain 
\begin{equation*}
\begin{split}
W(t) & \ge Cz_N^{-1}\mu(B(z,\sigma))
 +C\int_{\sigma^2}^t s^{-\frac{N(p-1)}{2}}W(s)^p\,ds
\end{split}
\end{equation*}
for a.a.~$t\in(\sigma^2,T/32)$. 

On the other hand, 
for a.a.~$z\in\Omega_{\sqrt{T}}$ and a.a.~$t\in(0,T/2)$, 
we observe from Lemma~\ref{Lemma:2.4} that $W(t)<\infty$. 
Applying Lemma~\ref{Lemma:2.6}, we obtain 
\[
	z_N^{-1}\mu(B(z,\sigma))\le C
	(\sigma^2)^{\frac{N}{2}-\frac{1}{p-1}}
	=C\sigma^{N-\frac{2}{p-1}}
\]
for all $\sigma\in(0,\sqrt{T}/16)$ and a.a.~$z\in\Omega_{\sqrt{T}}$, 
so that \eqref{eq:3.1} holds 
for all $\sigma\in(0,\sqrt{T}/16)$ 
and $z\in\Omega_{\sqrt{T}}$. 
Furthermore, in the case of $p=p_N$, 
we have  
\[
z_N^{-1}\mu(B(z,\sigma))\le C
\left[\log\frac{T}{2\sigma^2}\right]^{-\frac{N}{2}}
\le C
\left[\log\left(e+\frac{\sqrt{T}}{\sigma}\right)\right]^{-\frac{N}{2}}
\]
for all $\sigma\in(0,\sqrt{T}/16)$ 
and a.a.~$z\in\Omega_{\sqrt{T}}$. 
This implies that \eqref{eq:3.2} holds 
for all $\sigma\in(0,\sqrt{T}/16)$ 
and $z\in\Omega_{\sqrt{T}}$. 
Thus Proposition~\ref{Proposition:3.1} follows. 
$\Box$\vspace{5pt}

Next, we prove the following proposition.
\begin{proposition}
\label{Proposition:3.2}
Assume that there exists a supersolution to problem~\eqref{eq:P}
in $Q_T$, where $T>0$. 
Then there exist $C>0$ and $\epsilon\in(0,1)$ such that 
\begin{equation*}
\mu(B'(z',\sigma)\times[0,\sigma))
\le C\sigma^{N+1-\frac{2}{p-1}}
\end{equation*}
for all $z'\in{\mathbb R}^{N-1}$ and $\sigma\in(0,\epsilon\sqrt{T})$. 
\end{proposition}
For the proof of Proposition~\ref{Proposition:3.2}, 
we prepare the following lemma.
\begin{lemma}
\label{Lemma:3.3}
Let $u$ be a solution to problem~\eqref{eq:P} in $Q_T$, where $T>0$. 
Then there exists $C>0$ such that
\begin{equation*}
u(x,(2\sigma)^2) \geq C \sigma^{-N-1}
\mu(B'(z', \sigma)\times[0,\sigma))
\end{equation*}
for all $z'\in{\mathbb R}^{N-1}$, 
a.a.~$x\in B'(z',\sigma)\times(2\sigma,4\sigma)$, 
and a.a.~$\sigma\in(0,\sqrt{T}/16)$. 
\end{lemma}
{\bf Proof.} 
For any $x=(x',x_N)\in B'(z',\sigma)\times(2\sigma,4\sigma)$ and $y=(y',y_N)\in B'(z',\sigma)\times(0,\sigma)$, 
by \eqref{eq:1.4} and \eqref{eq:1.6} we apply the mean value theorem 
to find $\tilde{y}_N\in(0,y_N)$ so that
\begin{equation*}
\begin{split}
K(x,y,(2\sigma)^2) & =
y_N^{-1}G(x,y,(2\sigma)^2)=y_N^{-1}\Gamma_N(x-y,(2\sigma)^2)
 \left(1-\exp\left(-\frac{x_Ny_N}{(2\sigma)^2}\right)\right)\\
  & =\Gamma_N(x-y,(2\sigma)^2)\frac{x_N}{(2\sigma)^2}\exp\left(-\frac{x_N\tilde{y}_N}{(2\sigma)^2}\right)
  \ge C\sigma^{-N-1}.
\end{split}
\end{equation*}
Furthermore, 
by \eqref{eq:1.5} and \eqref{eq:1.6} we have
\begin{equation*}
\begin{split}
K(x,y,(2\sigma)^2) & =\frac{x_N}{(2\sigma)^2} \Gamma_N(x-y,(2\sigma)^2)
\ge C\sigma^{-N-1}
\end{split}
\end{equation*}
for $x\in B'(z',\sigma)\times(2\sigma,4\sigma)$ and $y\in B'(z',\sigma)\times\{0\}$. 
Then it follows from Definition~\ref{Definition:1.1} that
$$
u(x,(2\sigma)^2)\geq \int_{B'(z',\sigma)\times[0,\sigma)} K(x,y,(2\sigma)^2) \,d\mu(y)
\ge C\sigma^{-N-1}\mu(B'(z', \sigma)\times[0,\sigma))
$$
for all $z'\in{\mathbb R}^{N-1}$, a.a.~$x\in B'(z',\sigma)\times(2\sigma,4\sigma)$, and a.a.~$\sigma\in(0,\sqrt{T}/16)$. 
Thus Lemma~\ref{Lemma:3.3} follows. 
$\Box$\vspace{5pt}
\newline
{\bf Proof of Proposition~\ref{Proposition:3.2}.}
Assume that there exists a supersolution to problem~\eqref{eq:P}
in $Q_T$, where $T>0$. 
By Lemma~\ref{Lemma:2.5} we find a solution~$u$ to problem~\eqref{eq:P} in $Q_T$.

Let $\epsilon\in(0,1/16)$. 
For $\sigma\in(0,\epsilon\sqrt{T})$, we have 
\[
T-(2\sigma)^2>(1-4\epsilon^2)T>\frac{T}{2}.
\]
Set $\tilde{u}(x,t):=u(x,t+(2\sigma)^2)$. 
Then, by Lemma \ref{Lemma:2.4}, 
for a.a.~$\sigma\in(0,\epsilon\sqrt{T})$, the function 
$\tilde{u}$ is a solution to problem~\eqref{eq:P} 
with $\mu=x_Nu(x,(2\sigma)^2)$ in $Q_{T/2}$. 

Let $z'\in{\mathbb R}^{N-1}$ and set $z:=(z',3\sigma)\in \Omega_{3\sigma}$. 
Let $\delta\in(0,3/16)$. 
Since $\epsilon \sqrt{T}<\sqrt{T}/16$ and 
$B(z,\sigma)\subset B'(z',\sigma)\times(2\sigma,4\sigma)$, by Lemma~\ref{Lemma:3.3} we have 
\begin{equation}
\label{eq:3.5}
(3\sigma)^{-1}\int_{B(z,\delta\sigma)}y_Nu(y,(2\sigma)^2)\,dy
\ge C\sigma^{-1}\mu(B'(z', \sigma)\times[0,\sigma)).
\end{equation}
On the other hand, applying Proposition~\ref{Proposition:3.1} with $T=9\sigma^2$ to $\tilde{u}$, 
we have 
\[
(3\sigma)^{-1}\int_{B(z,\delta\sigma)}y_Nu(y,(2\sigma)^2)\,dy
=(3\sigma)^{-1}\int_{B(z,\delta\sigma)}y_N\tilde{u}(y,0)\,dy
\le C\sigma^{N-\frac{2}{p-1}}.
\]
This together with \eqref{eq:3.5} implies that
\begin{equation*}
\mu(B'(z', \sigma)\times[0,\sigma)) 
\le C\sigma^{N+1-\frac{2}{p-1}}
\end{equation*}
for all $z'\in{\mathbb R}^{N-1}$ and a.a.~$\sigma\in(0,\epsilon\sqrt{T})$. 
Then we obtain the desired inequality for all $z'\in{\mathbb R}^{N-1}$ and all $\sigma\in(0,\epsilon\sqrt{T})$. 
Thus Proposition~\ref{Proposition:3.2} follows. 
$\Box$\vspace{5pt}

Combining Propositions~\ref{Proposition:3.1} and \ref{Proposition:3.2}, 
we have: 
\begin{proposition}
\label{Proposition:3.3}
Assume that there exists a supersolution to problem~\eqref{eq:P}
in $Q_T$, where $T>0$. 
\begin{itemize}
  \item[{\rm (i)}]
  There exists $\gamma_1=\gamma_1(N,p)>0$ such that
  \[
  \mu(B_\Omega(z,\sigma)) \leq \gamma_1
  \sigma^{-\frac{2}{p-1}}\int_{B_\Omega(z,\sigma)}y_N\,dy
  \]
  for all $z\in\overline{\Omega}$ and $\sigma\in(0,\sqrt{T})$.  
  \item[{\rm (ii)}] 
  Let $p=p_N$. Then there exists $\gamma_2=\gamma_2(N)>0$ such that
  \[
  z_N^{-1}\mu(B(z,\sigma))
  \le \gamma_2
  \left[\log\left(e+\frac{\sqrt{T}}{\sigma}\right)\right]^{-\frac{N}{2}}
  \]
  for all $z=(z',z_N)\in\Omega_{3\sigma}$ and $\sigma\in(0,\sqrt{T})$. 
\end{itemize}
\end{proposition}
{\bf Proof.}
By Propositions~\ref{Proposition:3.1} and \ref{Proposition:3.2} 
we find $\delta\in(0,1/3)$ such that  
\begin{equation}
\label{eq:3.6}
\sup_{z\in\Omega_{\sigma}}\,z_N^{-1}\mu(B(z,\delta\sigma))
\le C \sigma^{N-\frac{2}{p-1}},\quad
\sup_{z\in\partial\Omega}\mu(B_\Omega(z,\delta\sigma))
\le C \sigma^{N+1-\frac{2}{p-1}},
\end{equation}
for all $\sigma\in(0,\sqrt{T})$. 
Furthermore, if $p=p_N$, then 
\begin{equation}
\label{eq:3.7}
\sup_{z\in\Omega_{\sigma}}\,z_N^{-1}\mu(B(z,\delta\sigma))
\le C \left[\log\left(e+\frac{\sqrt{T}}{\sigma}\right)\right]^{-\frac{N}{2}}
\quad\mbox{for all $\sigma\in(0,\sqrt{T})$}.
\end{equation}

Let $\sigma\in(0,\sqrt{T})$, $z=(z',z_N)\in\overline{\Omega}$, and $\overline{z}:=(z',0)\in\partial\Omega$.
Consider the case of $0\le z_N\le \delta\sigma/2$.  
Since $0<\delta<1/3$, 
we have
\[
B_\Omega(z,\delta\sigma/2)\subset B_\Omega(\overline{z},\delta\sigma)\subset B_\Omega(z,\sigma).
\]
Then, by \eqref{eq:3.6} we obtain
\begin{equation}
\label{eq:3.8}
\begin{split}
 \mu(B_\Omega(z,\delta\sigma/2))
  & \le \mu(B_\Omega(\overline{z},\delta\sigma))
 \le C\sigma^{N+1-\frac{2}{p-1}}\\
  & \le C\sigma^{-\frac{2}{p-1}}
 \int_{B_\Omega(\overline{z},\delta\sigma)}y_N\,dy
 \le C \sigma^{-\frac{2}{p-1}}\int_{B_\Omega(z,\sigma)}y_N\,dy.
\end{split}
\end{equation}
Consider the case of $z_N>\delta\sigma/2$. 
Then, by \eqref{eq:3.6} we have 
\begin{equation}
\label{eq:3.9}
\begin{split}
\mu(B_\Omega(z,\delta^2\sigma))
 & \le Cz_N\sigma^{N-\frac{2}{p-1}}
 \le Cz_N\sigma^{-\frac{2}{p-1}}\int_{B_\Omega(z,\delta^2\sigma/4)}\,dy\\
 & \le C\sigma^{-\frac{2}{p-1}}\int_{B_\Omega(z,\delta^2\sigma/4)}y_N\,dy
\le C\sigma^{-\frac{2}{p-1}}\int_{B_\Omega(z,\sigma)}y_N\,dy.
\end{split}
\end{equation}
Combining \eqref{eq:3.8} and \eqref{eq:3.9},
we obtain 
\begin{equation}
\label{eq:3.10}
\mu(B_\Omega(z,\delta^2\sigma/2))
\le C\sigma^{-\frac{2}{p-1}}\int_{B_\Omega(z,\sigma)}y_N\,dy
\end{equation}
for $z\in\overline{\Omega}$ and $\sigma\in(0,\sqrt{T})$. 
Therefore, by Lemma~\ref{Lemma:2.1}~(ii) and \eqref{eq:3.10}, 
for any $z\in\overline{\Omega}$, 
we find $\{\overline{z}_i\}_{i=1}^{m'}\subset B_\Omega(z,2\sigma)$ such that
\begin{equation*}
\begin{split}
\mu(B_\Omega(z,\sigma))
 & \le\sum_{i=1}^{m'}\mu(B_\Omega(\overline{z}_i,\delta^2\sigma/2))
\le C\sigma^{-\frac{2}{p-1}}\sum_{i=1}^{m'}\int_{B_\Omega(\overline{z}_i,\sigma)}y_N\,dy\\
 & \le C\sigma^{-\frac{2}{p-1}}\int_{B_\Omega(z,3\sigma)}y_N\,dy
 \le C\sigma^{-\frac{2}{p-1}}\int_{B_\Omega(z,\sigma)}y_N\,dy.
\end{split}
\end{equation*}
This implies assertion~(i).

Similarly, if $p=p_N$, then, by Lemma~\ref{Lemma:2.1}, 
for any $z\in\Omega_{3\sigma}$, 
we find $\{\tilde{z}_i\}_{i=1}^{m'}\subset B_\Omega(z,2\sigma)$ such that
$$
\mu(B(z,\sigma))\le\sum_{i=1}^{m'} \mu(B(\tilde{z}_i,\delta\sigma)).
$$
Since $\tilde{z}_i\in\Omega_\sigma$ and $0<\delta<1/3$, 
we deduce from \eqref{eq:3.7} that 
$$
z_N^{-1}\mu(B(z,\sigma))\le C\sum_{i=1}^{m'}\frac{z_N+2\sigma}{z_N}\left[\log\left(e+\frac{\sqrt{T}}{\sigma}\right)\right]^{-\frac{N}{2}}
\le C\left[\log\left(e+\frac{\sqrt{T}}{\sigma}\right)\right]^{-\frac{N}{2}}
$$
for all $z\in\Omega_{3\sigma}$ and $\sigma\in(0,\sqrt{T})$. 
This implies assertion~(ii). 
Thus Proposition~\ref{Proposition:3.3} follows.
$\Box$\vspace{5pt}

Next, we prove the following proposition.
\begin{proposition}
\label{Proposition:3.4}
Let $p=p_{N+1}$. 
Assume that there exists a supersolution to problem~\eqref{eq:P}
in~$Q_T$, where $T>0$. 
Then there exists $\gamma=\gamma(N)>0$ such that 
\begin{equation*}
\mu(B_\Omega(z,\sigma))\le \gamma
\left[\log\left(e+\frac{\sqrt{T}}{\sigma}\right)\right]^{-\frac{N+1}{2}}
\end{equation*}
for all $z\in\partial\Omega$ and $\sigma\in(0,\sqrt{T})$.
\end{proposition}
{\bf Proof.}
Let $p=p_{N+1}$. 
Assume that there exists a supersolution to problem~\eqref{eq:P}
in $Q_T$, where $T>0$. 
By Lemma~\ref{Lemma:2.5} we find a solution~$u$ to problem~\eqref{eq:P} in $Q_T$. 

Let $z=(z',0)\in\partial\Omega$.  
By Lemma \ref{Lemma:2.4}, 
for a.a.~$\sigma\in(0,\sqrt{T}/3)$, 
the function $v(x,t):= u(x,t+(2\sigma)^2)$
is a solution to problem~\eqref{eq:P} in $Q_{T-(2\sigma)^2}$. 
It follows from Proposition~\ref{Proposition:3.3} that 
\begin{equation}
\label{eq:3.11}
\int_{B_\Omega(z,r)} y_Nv(y,t)\,dy\le C
r^{-\frac{2}{p-1}}\int_{B_\Omega(z,r)}y_N\,dy
\end{equation}
for all $r\in(0,\sqrt{T-(2\sigma^2)-t})$ and a.a.~$t\in(0,T-(2\sigma)^2)$. 
Then 
$$
V(t):=t^{\frac{N}{2}+1}\int_\Omega  K(x,z,t)v(x,t)\,dx <\infty
\quad\mbox{for a.a.~$\displaystyle{t\in\left(\sigma^2,\frac{T-(2\sigma)^2}{2}\right)}$}.
$$
Indeed, 
by Lemma~\ref{Lemma:2.2} and \eqref{eq:3.11} we have
\begin{align*}
 \int_\Omega K(x,z,t) v(x,t)\,dx
  & \le Ct^{-1}\int_\Omega \Gamma_N(x-z,2t)x_Nv(x,t)\,dx\\
 & \le Ct^{-\frac{N}{2}-1}\sup_{z\in\overline{\Omega}}
 \int_{B_\Omega(z,\sqrt{2t})}y_N v(y,t)\,dy<\infty
\end{align*}
for a.a.~$t\in(\sigma^2,(T-(2\sigma)^2)/2)$. 

We derive an integral inequality for $V$. 
By Fubini's theorem, \eqref{eq:1.5}, and \eqref{eq:2.6} we have
\begin{equation}
\label{eq:3.12}
\begin{split}
  & \int_\Omega K(x,z,t)v(x,t)\,dx\\
  & \ge\int_\Omega\int_\Omega K(x,z,t)G(x,y,t)v(y,0)\,dx\,dy\\
  & \qquad\quad
  +\int_0^t\int_\Omega\int_\Omega K(x,z,t)G(x,y,t-s)v(y,s)^p\,dx\,dy\,ds\\
  & =\int_\Omega \frac{y_N}{2t}\Gamma_N(y-z,2t)v(y,0)\,dy
  +\int_0^t \int_\Omega 
\frac{y_N}{2t-s}\Gamma_N(y-z,2t-s) v(y,s)^p \,dy\,ds.
\end{split}
\end{equation}
Furthermore, 
\begin{equation}
\label{eq:3.13}
\begin{split}
\int_\Omega \frac{y_N}{2t}\Gamma_N(y-z,2t)v(y,0)\,dy
 & \ge\int_{B'(z',\sigma)\times(2\sigma,4\sigma)} \frac{y_N}{2t}
 \Gamma_N(y-z,2t)u(y,(2\sigma)^2)\,dy\\
 & \ge Ct^{-\frac{N}{2}-1}\int_{B'(z',\sigma)\times(2\sigma,4\sigma)} y_Nu(y,(2\sigma)^2)\,dy
\end{split}
\end{equation}
for all $t\in(\sigma^2,(T-(2\sigma)^2)/3)$. 
On the other hand, 
by Lemma~\ref{Lemma:3.2}~(i) and \eqref{eq:1.5} 
we obtain
\begin{equation*}
\begin{split}
 & \int_0^t \int_\Omega 
\frac{y_N}{2t-s}\Gamma_N(y-z,2t-s) v(y,s)^p \,dy\,ds\\
 & \ge\int_0^t \int_\Omega 
\frac{y_N}{2t}\left(\frac{s}{2t}\right)^{\frac{N}{2}}\Gamma_N(y-z,s) v(y,s)^p \,dy\,ds
=\int_0^t \int_\Omega 
\left(\frac{s}{2t}\right)^{\frac{N}{2}+1}K(y,z,s) v(y,s)^p \,dy\,ds.
\end{split}
\end{equation*}
Then Jensen's inequality together with \eqref{eq:2.5} implies that
\begin{equation}
\label{eq:3.14}
\begin{split}
 & \int_0^t \int_\Omega 
\frac{y_N}{2t-s}\Gamma_N(y-z,2t-s) v(y,s)^p \,dy\,ds\\
 & \ge\int_0^t  
\left(\frac{s}{2t}\right)^{\frac{N}{2}+1}(\pi s)^{-\frac{1}{2}}
\left(\int_\Omega (\pi s)^{\frac{1}{2}}K(y,z,s) v(y,s) \,dy\right)^p\,ds\\
 & \ge Ct^{-\frac{N}{2}-1}\int_0^t s^{\frac{N+1}{2}}
 \left(\int_\Omega \sqrt{s}K(y,z,s) v(y,s) \,dy\right)^p\,ds.
\end{split}
\end{equation}
Since $p=p_{N+1}=1+2/(N+1)$, 
by \eqref{eq:3.12}, \eqref{eq:3.13}, and \eqref{eq:3.14}
we see that
\begin{equation}
\label{eq:3.15}
\begin{split}
V(t) & \ge C\int_{B'(z',\sigma)\times(2\sigma,4\sigma)} y_Nu(y,(2\sigma)^2)\,dy
+C\int_0^t s^{-1}V(s)^p\,ds\\
 & \ge C\int_{B'(z',\sigma)\times(2\sigma,4\sigma)} y_Nu(y,(2\sigma)^2)\,dy
+C\int_{t_*}^t s^{-1}V(s)^p\,ds
\end{split}
\end{equation}
for a.a.~$t\in(\sigma^2,(T-(2\sigma)^2)/3)$, all $t_*\in(0,t)$, and a.a.~$\sigma\in(0,\sqrt{T}/3)$.

Let $\epsilon>0$ be small enough. 
We apply Lemma~\ref{Lemma:2.6} with $t_*=t/2$ to inequality~\eqref{eq:3.15}. 
Then 
\begin{equation}
\label{eq:3.16}
\int_{B'(z',\sigma)\times(2\sigma,4\sigma)}  y_Nu(y,(2\sigma)^2) \,dy 
\leq C \left[\log\left(e+\frac{\sqrt{T}}{\sigma}\right)\right]^{-\frac{N+1}{2}}
\end{equation}
for a.a.~$\sigma\in(0,\epsilon \sqrt{T})$.
Furthermore, by Lemma~\ref{Lemma:3.3}, 
taking small enough~$\epsilon>0$ if necessary, we have
\begin{equation}
\label{eq:3.17}
\begin{split}
&\int_{B'(z',\sigma)\times(2\sigma,4\sigma)}  y_Nu(y,(2\sigma)^2) \,dy\\
&\geq C \sigma^{-N-1} \mu(B'(z', \sigma)\times[0,\sigma))\int_{B'(z',\sigma)\times(2\sigma,4\sigma)} y_N\,dy
\geq C\mu(B'(z', \sigma)\times[0,\sigma))
\end{split}
\end{equation}
for a.a.~$\sigma\in(0,\epsilon\sqrt{T})$. 
Combining \eqref{eq:3.16} and \eqref{eq:3.17}, 
we find $\delta\in(0,1)$ such that 
\[
\sup_{z\in\partial\Omega}\mu(B_\Omega(z, \delta\sigma))
\le C\left[\log\left(e+\frac{\sqrt{T}}{\sigma}\right)\right]^{-\frac{N+1}{2}}
\quad\mbox{for a.a.~$\sigma\in(0,\sqrt{T})$}.
\]
This together with Lemma~\ref{Lemma:2.1} implies that
$$
\sup_{z\in\partial\Omega}\mu(B_\Omega(z, \sigma))
\le\sum_{i=1}^{m'}\mu(B_\Omega(z_i', \delta\sigma)) 
\le C\left[\log\left(e+\frac{\sqrt{T}}{\sigma}\right)\right]^{-\frac{N+1}{2}}
$$
for all~$\sigma\in(0,\sqrt{T})$. 
Thus Proposition~\ref{Proposition:3.4} follows.
$\Box$\vspace{5pt}

Now we are ready to complete the proof of Theorem~\ref{Theorem:1.2}.
\vspace{5pt}
\newline
{\bf Proof of Theorem~\ref{Theorem:1.2}.} 
Let $u$ be a supersolution to problem~\eqref{eq:P} in $Q_T$, where $T>0$.
By Propositions~\ref{Proposition:3.3} and \ref{Proposition:3.4} 
we have only to show that $\mu(\partial\Omega)=0$ if $p\geq2$. 
Furthermore, thanks to Remark~\ref{Remark:1.2}, 
it suffices to consider the case of $p=2$.

Let $p=2$. Consider the case of $N=1$. 
Proposition~\ref{Proposition:3.4} implies that
$$
\mu(\partial\Omega)=\mu(\{0\})\leq C
\left[\log\left(e+\frac{\sqrt{T}}{\sigma}\right)\right]^{-1}\to 0
\quad\mbox{as }\sigma\to 0.
$$
This means that $\mu(\partial\Omega)=0$. 
Thus Theorem~\ref{Theorem:1.2} follows in the case of $N=1$.

Consider the case of $N\geq2$. Let $z=(z',0)\in\partial\Omega$. 
Set
\[
U(x_N,t):= \int_{{\mathbb R}^{N-1}} \Gamma_{N-1}(x',T')u(x+z,t) \,dx',
\quad\mbox{where}\quad T'=\frac{T}{4}.
\]
By Lemma~\ref{Lemma:2.4}, \eqref{eq:1.3}, and \eqref{eq:1.4}, 
for any $L_1$, $L_2>0$ with $L_1<L_2$, 
we have 
\begin{align*}
\int_{L_1}^{L_2} U(x_N,t)\,dx_N & 
\le CT^{\frac{N-1}{2}}
\int_{{\mathbb R}^{N-1}\times(L_1,L_2)}
G(y,x+z,2T')u(x+z,t)\,dx\\
 & \le CT^{\frac{N-1}{2}}\int_\Omega G(y,x,T')u(x,t)\,dx<\infty
\end{align*}
for a.a.~$y=(y',y_N)\in B'(0,T')\times(L_1,L_2)$ and a.a.~$t\in(0,T')$. 
This implies that 
$U(x_N,t)<\infty$ for a.a.~$(x_N,t)\in(0,\infty)\times(0,T')$.
Furthermore, it follows from Definition~\ref{Definition:1.1}  that
\begin{equation}
\label{eq:3.18}
\begin{split}
U(x_N,t) 
& \ge \int_{{\mathbb R}^{N-1}} \Gamma_{N-1}(x',T')
\left(\int_{\partial\Omega}K(x+z,y,t)\,d\mu(y')\right)\,dx'\\
& +\int_{{\mathbb R}^{N-1}}\int_0^t\int_\Omega\,\Gamma_{N-1}(x',T')G(x+z,y,t-s)u(y,s)^p\,dy\,ds\,dx'
\end{split}
\end{equation}
for a.a.~$(x_N,t)\in(0,\infty)\times(0,T')$. 
By \eqref{eq:1.5} we see that 
\begin{equation}
\label{eq:3.19}
\begin{split}
 & \int_{{\mathbb R}^{N-1}} \Gamma_{N-1}(x',T')
\left(\int_{\partial\Omega}K(x+z,y,t)\,d\mu(y')\right)\,dx'\\
 & =\frac{x_N}{t}\int_{{\mathbb R}^{N-1}}\int_{{\mathbb R}^{N-1}} \Gamma_{N-1}(x',T')\Gamma_N(x+z-y,t)\,dx'\,d\mu(y')\\
 & =\frac{x_N}{t}\Gamma_1(x_N,t)\int_{{\mathbb R}^{N-1}}\Gamma_{N-1}(z'-y',T'+t)\,d\mu(y')\\
 & \ge \frac{x_N}{t}\Gamma_1(x_N,t)\int_{B'(z',\sqrt{T'})}\Gamma_{N-1}(z'-y',T'+t)\,d\mu(y')\\
 & \ge C(T')^{-\frac{N-1}{2}}\mu(B'(z',\sqrt{T'}))\frac{x_N}{t}\Gamma_1(x_N,t)
\end{split}
\end{equation}
for all $(x_N,t)\in(0,\infty)\times(0,T')$. 
On the other hand, 
by \eqref{eq:1.4} and Jensen's inequality we have 
\begin{equation}
\label{eq:3.20}
\begin{split}
 & \int_{{\mathbb R}^{N-1}}\int_0^t\int_\Omega\,\Gamma_{N-1}(x',T')G(x+z,y,t-s)u(y,s)^p\,dy\,ds\,dx'\\
 & =\int_{{\mathbb R}^{N-1}}\int_0^t\int_0^\infty\int_{{\mathbb R}^{N-1}}\Gamma_{N-1}(x',T')
\Gamma_{N-1}(x'+z'-y',t-s)\\
 &\quad\times
 \left[\Gamma_1(x_N-y_N,t-s)- \Gamma_1(x_N+y_N,t-s)\right]u(y,s)^p\,dy'\,dy_N\,ds\,dx'\\
 & =\int_0^t\int_0^\infty\int_{{\mathbb R}^{N-1}}\Gamma_{N-1}(y'-z',T'+t-s)u(y,s)^p\\
 &\quad\times
 \left[\Gamma_1(x_N-y_N,t-s)- \Gamma_1(x_N+y_N,t-s)\right]\,dy'\,dy_N\,ds\\
  & =\int_0^t\int_0^\infty\int_{{\mathbb R}^{N-1}}\Gamma_{N-1}(y',T'+t-s)u(y+z,s)^p\\
 &\quad\times
 \left[\Gamma_1(x_N-y_N,t-s)- \Gamma_1(x_N+y_N,t-s)\right]\,dy'\,dy_N\,ds
\end{split}
\end{equation}
for all $(x_N,t)\in(0,\infty)\times(0,T')$. 
Since 
\[
\frac{1}{2}\leq \frac{T'}{T'+t-s}\leq 1\quad\mbox{for $t$, $s\in(0,T')$ with $s<t$},
\]
we see that 
\[
\Gamma_{N-1}(y',T'+t-s) \geq \left(\frac{T'}{T'+t-s}\right)^\frac{N-1}{2} 
\Gamma_{N-1}(y',T') \geq C \Gamma_{N-1}(y',T')
\]
for all $y'\in {\mathbb R}^{N-1}$ and $t$, $s\in(0,T')$ with $s<t$.
This together with \eqref{eq:3.20} and Jensen's inequality implies that 
\begin{equation}
\label{eq:3.21}
\begin{split}
&\int_{{\mathbb R}^{N-1}}\int_0^t\int_\Omega\,\Gamma_{N-1}(x',T')G(x+z,y,t-s)u(y,s)^p\,dy\,ds\,dx' \\
&\geq
C \int_0^t\int_0^\infty \left[\Gamma_1(x_N-y_N,t-s)- \Gamma_1(x_N+y_N,t-s)\right]U(y_N,s)^p\,dy_N\,ds
\end{split}
\end{equation}
for all $(x_N,t)\in(0,\infty)\times(0,T')$. 
Combining \eqref{eq:3.18}, \eqref{eq:3.19}, and \eqref{eq:3.21}, 
we obtain 
\begin{equation*}
\begin{split}
U(x_N,t) &\geq C(T')^{-\frac{N-1}{2}} \mu(B'(z',\sqrt{T'}))  \frac{x_N}{t}\Gamma_1 (x_N,t) \\
&+ C\int_0^t \int_0^\infty [\Gamma_1(x_N-y_N,t-s)-\Gamma_1(x_N+y_N,t-s)] U(y_N,s)^p\,dy_Nds
\end{split}
\end{equation*}
for a.a.~$(x_N,t)\in(0,\infty)\times(0,T')$. 
This means that $U$ is a supersolution to problem~\eqref{eq:P} in~$Q_{T'}$ 
with $N=1$ and the initial data $C(T')^{-\frac{N-1}{2}} \mu(B'(z',\sqrt{T'}))\delta_1$. 
Then we apply Theorem~\ref{Theorem:1.2}~(iii) with $N=1$ to obtain $\mu(B'(z',\sqrt{T'}))=0$. 
Since $z'\in{\mathbb R}^{N-1}$ is arbitrary, we deduce that $\mu(\partial\Omega)=0$. 
Therefore the proof of Theorem~\ref{Theorem:1.2} is complete. 
$\Box$
\section{Proof of Theorem~\ref{Theorem:1.1}}
%
For the proof of Theorem~\ref{Theorem:1.1}, 
we prepare the following lemma.
\begin{lemma}
\label{Lemma:4.1} 
Let $\phi\in C_0^\infty({\mathbb R}^N)$. Set 
$$
\varphi_0(x):=x_N\phi(x),\quad
\varphi(x,t):=\int_\Omega G(x,y,t)\varphi_0(y)\,dy,
\quad\mbox{for $(x,t)\in Q_\infty$}.
$$
Then
\begin{equation}
\label{eq:4.1}
\lim_{t\to +0}\sup_{x\in\Omega}|x_N^{-1}\varphi(x,t)-\phi(x)|=0,
\quad
\lim_{t\to +0}\sup_{x\in\partial\Omega}|(\partial_{x_N}\varphi)(x,t)-\phi(x)|=0.
\end{equation}
Furthermore, for any $T>0$, there exists $C>0$ such that 
\begin{equation}
\label{eq:4.2}
\sup_{(x,t)\in Q_T}|x_N^{-1}\varphi(x,t)|
+\sup_{(x,t)\in\partial\Omega\times(0,T)}|(\partial_{x_N}\varphi)(x,t)|\le C\exp\left(-\frac{|x|^2}{CT}\right)
\quad\mbox{for $(x,t) \in Q_T$}.
\end{equation}
\end{lemma}
{\bf Proof.}
Since $\phi\in C_0^\infty({\mathbb R}^N)$ and $\varphi_0=0$ on $\partial\Omega$, 
by parabolic regularity theorems we see that $\varphi\in C^{2,1}(\overline{Q_\infty})$, which implies that 
\begin{equation}
\label{eq:4.3}
\lim_{t\to +0}\sup_{x\in B_\Omega(0,L)}|(\partial_{x_N}\varphi)(x,t)-(\partial_{x_N}\varphi_0)(x)|=0\quad\mbox{for all $L>0$}. 
\end{equation}
For any $x=(x',x_N)\in \Omega$ and $t>0$,
we have 
\begin{align}
\label{eq:4.4}
x_N^{-1}\varphi(x,t) & =x_N^{-1}\int_0^1 \frac{d}{ds}\varphi(x',sx_N,t)\,ds
=\int_0^1 (\partial_{x_N}\varphi)(x',sx_N,t)\,ds,\\
\nonumber
\phi(x)=x_N^{-1}\varphi_0(x) & =x_N^{-1}\int_0^1 \frac{d}{ds}\varphi_0(x',sx_N)\,ds
=\int_0^1 (\partial_{x_N}\varphi_0)(x',sx_N)\,ds.
\end{align}
These together with \eqref{eq:4.3} imply that 
\begin{equation}
\label{eq:4.5}
\lim_{t\to +0}\sup_{x\in \Omega\cap B(0,L)}\left|x_N^{-1}\varphi(x,t)-\phi(x)\right|=0\quad\mbox{for all $L>0$}.
\end{equation}

Let $R>0$ be such that $\mbox{supp}\,\phi\subset B(0,R)$. 
Since $|x-y|\ge |x|/2$ for $x\in\Omega_{2R}$ and $y\in B_\Omega(0,R)$, 
it follows from \eqref{eq:1.4} that 
\begin{equation}
\label{eq:4.6}
\begin{split}
|x_N^{-1}\varphi(x,t)| & \le Ct^{-\frac{N}{2}}\int_{B_\Omega(0,R)}\exp\left(-\frac{|x-y|^2}{4t}\right)\,dy\\
 & \le Ct^{-\frac{N}{2}}\exp\left(-\frac{|x|^2}{Ct}\right)\int_{{\mathbb R}^N}\exp\left(-\frac{|x-y|^2}{8t}\right)\,dy
\le C\exp\left(-\frac{|x|^2}{Ct}\right)
\end{split}
\end{equation}
for $x=(x',x_N)\in\Omega_{2R}$ and $t>0$. 
Similarly, by \eqref{eq:1.4} we have 
\begin{equation}
\label{eq:4.7}
\begin{split}
 & |(\partial_{x_N}\varphi)(x,t)|
=\left|\int_\Omega (\partial_{x_N}G)(x,y,t)\varphi_0(y)\,dy\right|\\
 & \le C\int_{B_\Omega(0,R)}
\left|(\partial_{x_N}\Gamma_N)(x-y,t)\left(1-\exp\left(-\frac{x_N y_N}{t}\right)\right)
+\Gamma_N(x-y,t)\frac{y_N}{t}\exp\left(-\frac{x_Ny_N}{t}\right)\right|\,dy\\
 & \le Ct^{-\frac{N}{2}}\int_{B_\Omega(0,R)} \left(\frac{|x-y|}{t}+\frac{|y|}{t}\right)\exp\left(-\frac{|x-y|^2}{4t}\right)\,dy\\
 & \le CR^{-1}t^{-\frac{N}{2}}\int_{B_\Omega(0,R)} \left(\frac{|x-y|^2}{t}+\frac{|x|^2}{t}\right)\exp\left(-\frac{|x-y|^2}{4t}\right)\,dy\\
 & \le Ct^{-\frac{N}{2}}\exp\left(-\frac{|x|^2}{Ct}\right)\int_{{\mathbb R}^N}\exp\left(-\frac{|x-y|^2}{8t}\right)\,dy\le C\exp\left(-\frac{|x|^2}{Ct}\right)
\end{split}
\end{equation}
for all $x\in\overline{\Omega}\setminus B_\Omega(0,2R)$ and $t>0$. 
Then, by \eqref{eq:4.4} and \eqref{eq:4.7}
we see that
\begin{equation}
\label{eq:4.8}
|x_N^{-1}\varphi(x,t)|\le C\exp\left(-\frac{|x'|^2}{Ct}\right)\le C\exp\left(-\frac{|x|^2}{Ct}\right)
\end{equation}
for all $(x',x_N,t)\in Q_\infty$ with $|x'|\ge 2R$ and $x_N\le 2R$. 
Since 
$$
\lim_{x_N\to +0}x_N^{-1}\varphi(x',x_N,t)=(\partial_{x_N}\varphi)(x',0,t)\quad\mbox{for $(x',0,t)\in\partial\Omega\times(0,\infty)$}, 
$$
by \eqref{eq:4.3}, \eqref{eq:4.5}, \eqref{eq:4.6}, and \eqref{eq:4.7} 
we obtain \eqref{eq:4.1}. 
Furthermore, thanks to $\varphi\in C^{2,1}(\overline{Q_\infty})$, 
by \eqref{eq:4.4} we have 
\begin{equation}
\label{eq:4.9}
|x_N^{-1}\varphi(x,t)|\le C\le C\exp\left(-\frac{|x|^2}{CT}\right)
\end{equation}
for all $(x',x_N,t)\in Q_T$ with $|x'|<2R$ and $x_N\le 2R$. 
Then, combining \eqref{eq:4.6}, \eqref{eq:4.8}, and \eqref{eq:4.9}, 
we deduce \eqref{eq:4.2}. The proof is complete.
$\Box$\vspace{5pt}

Next, we prove a lemma on the existence and the uniqueness of the initial trace of $x_Nu(x,t)$.
\begin{lemma}
\label{Lemma:4.2}
Let $u\in {\mathcal L}(Q_T)$ solve {\rm (E)} in $Q_T$, where $T>0$. 
Then 
\begin{equation}
\label{eq:4.10}
\underset{t\in(0,T-\epsilon)}{{\rm ess\,\,sup}}
\,\int_{B_\Omega(0,R)}y_Nu(y,t)\,dy<\infty
\end{equation}
for all $R>0$ and $\epsilon\in(0,T)$. 
Furthermore, 
there exists a unique $\mu\in{\mathcal M}$ such that
\begin{equation}
\label{eq:4.11}
\underset{t\to +0}{\mbox{{\rm ess lim}}}
\int_\Omega y_N u(y,t)\phi(y)\,dy=\int_{\overline{\Omega}} \phi(y)\,d\mu(y)
\quad\mbox{for all $\phi\in C_0({\mathbb R}^N)$}.
\end{equation}
\end{lemma}
{\bf Proof.}
Let $u\in {\mathcal L}(Q_T)$ solve {\rm (E)} in $Q_T$, where $T>0$. 
By Definition~\ref{Definition:1.1} we see that, 
for a.a.~$\tau\in(0,T)$, the function
$$
u_\tau(x,t):=u(x,t+\tau)\quad\mbox{for a.a.~$(x,t)\in Q_{T-\tau}$}, 
$$
 is a solution to problem~\eqref{eq:P} with $\mu=x_Nu(x,\tau)$ in $Q_{T-\tau}$. 
Then, by Theorem~\ref{Theorem:1.2}, 
for any $R>0$ and $\epsilon\in(0,T)$, we have
\begin{equation}
\label{eq:4.12}
\underset{t\in(0,T-\epsilon)}{{\rm ess\,\,sup}}
\int_{B_\Omega(z,\sigma)}y_Nu(y,\tau)\,dy
\leq C\sigma^{-\frac{2}{p-1}}\int_{B(z,\sigma)}y_N\,dy
\le C(|z|+\sqrt{T})\sigma^{N-\frac{2}{p-1}}
\end{equation}
for all $z\in \overline{\Omega}$ and $\sigma\in(0,\sqrt{\epsilon})$. 
This together with Lemma~\ref{Lemma:2.1} implies \eqref{eq:4.10}. 

We prove \eqref{eq:4.11}. 
By Definition~\ref{Definition:1.1}~(i) we find a measurable set $\Sigma_T\subset(0,T)$ such that 
\begin{itemize}
  \item 1-dimensional Lebesgue measure of the set $(0,T)\setminus\Sigma_T$ is zero;
  \item for any nonnegative function $\phi\in C_0({\mathbb R}^N)$ and any $t$, $s\in\Sigma_T$ with $t>s$,
  \begin{equation}
  \label{eq:4.13}
  \begin{split}
   \int_\Omega x_Nu(x,t)\phi(x)\,dx
   & \ge \int_\Omega\int_\Omega G(x,y,t-s)u(y,s)x_N\phi(x)\,dy\,dx\\
   & \qquad
   +\int_\Omega \int_s^t\int_\Omega 
   G(x,y,t-\tau)u(y,\tau)^p x_N\phi(x)\,dy\,d\tau\,dx. 
   \end{split}
  \end{equation}
\end{itemize}
Let $\{t_j\}$ be a sequence in $\Sigma_T$ such that $\lim_{j\to\infty}t_j=0$.
Applying the weak compactness of Radon measures (see e.g. \cite{EG}*{Section~1.9}), 
by \eqref{eq:4.12} 
we find a subsequence $\{t_{j_k}\}$ of $\{t_j\}$ and $\mu\in{\mathcal M}$ such that 
\begin{equation}
\label{eq:4.14}
\lim_{k\to\infty}\int_\Omega y_Nu(y,t_{j_k})\psi(y)\,dy=\int_{\overline{\Omega}}\psi(y)\,d\mu(y)
\end{equation}
for all $\psi\in C_0({\mathbb R}^N)$. 

Let $\{s_j\}$ be a sequence in $\Sigma_T$ such that $\lim_{j\to\infty}s_j=0$. 
Similarly to \eqref{eq:4.14}, 
we find a subsequence $\{s_{j_k}\}$ of $\{s_j\}$ and $\mu'\in{\mathcal M}$ such that 
\begin{equation}
\label{eq:4.15}
\lim_{k\to\infty}\int_{\Omega}y_Nu(y,s_{j_k})\psi(y)\,dy=\int_{\overline{\Omega}}\psi(y)\,d\mu'(y)
\end{equation}
for all $\psi\in C_0({\mathbb R}^N)$. 
Taking a subsequence if necessary, we can assume that $t_{j_k}>s_{j_k}$ for all $k=1,2,\dots$.
Let $\phi\in C_0^\infty({\mathbb R}^N)$ be such that $\phi\ge 0$ in ${\mathbb R}^N$ 
and $\mbox{supp}\,\phi\subset B(0,R)$ for some $R>0$. 
Let $\varphi_0$ and $\varphi$ be as in Lemma~\ref{Lemma:4.1}. 
It follows from \eqref{eq:4.13} that 
\begin{equation*}
\begin{split}
 & \int_\Omega x_Nu(x,t_{j_k})\phi(x)\,dx
 \ge\int_\Omega \left(\int_\Omega G(x,y,t_{j_k}-s_{j_k})\varphi_0(x)\,dx\right)\,u(y,s_{j_k})\,dy\\
 & =\int_\Omega \varphi(y,t_{j_k}-s_{j_k})u(y,s_{j_k})\,dy
 \ge\int_{B_\Omega(0,R)} \varphi(y,t_{j_k}-s_{j_k})u(y,s_{j_k})\,dy\\
 & \ge \int_\Omega y_N\phi(y)u(y,s_{j_k})\,dy-\sup_{x\in B_\Omega(0,R)}\left|x_N^{-1}\varphi(x,t_{j_k}-s_{j_k})-\phi(x)\right|
 \int_{B_\Omega(0,R)}y_Nu(y,s_{j_k})\,dy.
\end{split} 
\end{equation*}
Then, by \eqref{eq:4.1}, \eqref{eq:4.10}, \eqref{eq:4.14}, and \eqref{eq:4.15}
we obtain 
$$
\int_{\overline{\Omega}}\phi(y)\,d\mu(y)\ge\int_{\overline{\Omega}}\phi(y)\,d\mu'(y).
$$
Since $\phi$ is arbitrary, we deduce that $\mu\ge\mu'$ in ${\mathcal M}$. 
Similarly, we have  $\mu'\ge\mu$ in ${\mathcal M}$. Thus we see that $\mu=\mu'$ in ${\mathcal M}$. 
Since $\{s_j\}\subset\Sigma_T$ is arbitrary, we obtain \eqref{eq:4.11}.
Thus Lemma~\ref{Lemma:4.2} follows.
$\Box$\vspace{5pt}
\begin{lemma}
\label{Lemma:4.3} 
Let $u$ and $\mu$ be as in Lemma~{\rm\ref{Lemma:4.2}}. 
Then, for any $\epsilon\in(0,T)$ and $\delta>0$,
$$
\underset{t\in(0,T-\epsilon)}{{\rm ess\,\,sup}}\,\int_\Omega e^{-\delta|x|^2}x_Nu(x,t)\,dx+
\int_{\overline{\Omega}}e^{-\delta|x|^2}\,d\mu(x)<\infty. 
$$
\end{lemma}
{\bf Proof.}
Let $\epsilon\in(0,T)$ and set $\sigma=\sqrt{T}/2$. 
It follows from \eqref{eq:4.11} and \eqref{eq:4.12} that 
\begin{equation}
\label{eq:4.16}
\underset{t\in(0,T-\epsilon)}{{\rm ess\,\,sup}}\,\int_{B_\Omega(z,\sigma)} x_Nu(x,t)\,dx
+\mu(B_\Omega(z,\sigma))\le C\sigma^{N-\frac{2}{p-1}}(|z|+\sqrt{T}).
\end{equation}
By the Besicovitch covering lemma we find an integer $m=m(N)$ 
and $\{x_{k,i}\}_{k=1,\dots,m,\,i\in{\bf N}}\subset\overline{\Omega}$ such that 
\begin{equation}
\label{eq:4.17}
B_{k,i}\cap B_{k,j}=\emptyset\quad\mbox{if}\quad i\not=j,
\qquad
\overline{\Omega}\subset\bigcup_{k=1}^m\bigcup_{i=1}^\infty B_{k,i},
\end{equation}
where $B_{k,i}:=\overline{B(x_{k,i},\sigma)}$. 
By \eqref{eq:4.16} and \eqref{eq:4.17}, for any $\delta>0$, 
\begin{equation}
\label{eq:4.18}
\begin{split}
 & \int_{\overline{\Omega}} e^{-\delta|x|^2}x_Nu(x,t)\,dx
+\int_{\overline{\Omega}}e^{-\delta|x|^2}\,d\mu(x)\\
 & \le\sum_{k=1}^m\sum_{j=1}^\infty\left(\int_{B_{k,i}\cap\,\overline{\Omega}}
 e^{-\delta|x|^2}x_Nu(x,t)\,dx+\int_{B_{k,i}\cap\,\overline{\Omega}}e^{-\delta|x|^2}\,d\mu(x)\right)\\
 & \le\sum_{k=1}^m\sum_{j=1}^\infty\sup_{x\in B_{k,i}}e^{-\delta|x|^2}
\left(\int_{B_{k,i}\cap\,\overline{\Omega}} x_Nu(x,t)\,dx+\mu(B_{k,i}\cap\,\overline{\Omega})\right)\\
 & \le C\sigma^{N-\frac{2}{p-1}}\sum_{k=1}^m\sum_{j=1}^\infty\sup_{x\in B_{k,i}}e^{-\delta|x|^2}(|x_{k,i}|+\sqrt{T})
\end{split}
\end{equation}
for a.a.~$t\in(0,T-\epsilon)$. 
On the other hand, for any $x$, $y\in B_{k,i}$, we have
$$
e^{-\delta|x|^2}\le e^{-\delta ||y|-|x-y||^2}
=e^{-\delta\left(|y|^2+|x-y|^2-2|y||x-y|\right)}
\le Ce^{-\frac{\delta}{2}|y|^2},
$$
which implies that 
$$
\sup_{x\in B_{k,i}}e^{-\delta|x|^2}(|x_{k,i}|+\sqrt{T})
\le Ce^{-\frac{\delta}{2}|y|^2} 
(|x_{k,i}|+\sqrt{T})\\
\le Ce^{-\frac{\delta}{4}|y|^2}. 
$$
Then we obtain
\begin{equation}
\label{eq:4.19}
\sup_{x\in B_{k,i}}e^{-\delta|x|^2}(|x_{k,i}|+\sqrt{T})
\le C\sigma^{-N}\int_{B_{k,i}}e^{-\frac{\delta}{4}|x|^2}\,dx
\end{equation}
for $k=1,\dots,m$ and $i\in{\bf N}$. 
Therefore, by \eqref{eq:4.17}, \eqref{eq:4.18}, and \eqref{eq:4.19} we obtain 
\begin{align*}
 & \int_{\overline{\Omega}} e^{-\delta|x|^2}x_Nu(x,t)\,dx+\int_{\overline{\Omega}}e^{-\delta|x|^2}\,d\mu(x)\\
 & \le C\sigma^{-\frac{2}{p-1}}\sum_{k=1}^m\sum_{j=1}^\infty\int_{B_{k,i}}e^{-\frac{\delta}{4}|x|^2}\,dx
 \le C\sigma^{-\frac{2}{p-1}}\sum_{k=1}^m\int_{{\mathbb R}^N} e^{-\frac{\delta}{4}|x|^2}\,dx\le C\sigma^{-\frac{2}{p-1}}
\end{align*}
for a.a.~$t\in(0,T-\epsilon)$. Thus Lemma~\ref{Lemma:4.3}  follows. 
$\Box$
\vspace{5pt}

Now we are ready to complete the proof of Theorem~\ref{Theorem:1.1}.
We show assertion~(ii), and then we prove assertion~(i). 
\vspace{5pt}
\newline
{\bf Proof of Theorem~\ref{Theorem:1.1}~(ii).}
Let $\mu\in{\mathcal M}$ and let $u$ be a solution to problem~\eqref{eq:P} in $Q_T$, where $T>0$. 
Lemma~\ref{Lemma:2.4} implies that $u$ solves problem~(E) in $Q_T$.  
Then, by Lemma~\ref{Lemma:4.2} we find a unique $\nu\in{\mathcal M}$ such that 
\begin{equation}
\label{eq:4.20}
\underset{t\to +0}{\mbox{{\rm ess lim}}}
\int_\Omega y_N u(y,t)\psi(y)\,dy=\int_{\overline{\Omega}}\psi(y)\,d\nu(y),
\quad\mbox{for all $\psi\in C_0({\mathbb R}^N)$}.
\end{equation}

We prove that $\mu=\nu$ in ${\mathcal M}$. 
Let $\phi\in C_0^\infty({\mathbb R}^N)$ be such that $\phi\ge 0$ in ${\mathbb R}^N$ and $\mbox{supp}\,\phi\subset B(0,R)$ for some $R>0$. 
Let $\varphi_0$ and $\varphi$ be as in Lemma~\ref{Lemma:4.1}. 
Let $\Sigma_T$ be as in the proof of Lemma~\ref{Lemma:4.2}. 
Then \eqref{eq:4.13} holds. 
Furthermore, we have 
\begin{equation*}
\begin{split}
 & \int_\Omega\int_\Omega G(x,y,t-\tau)u(y,\tau)x_N\phi(x)\,dx\,dy
\ge \int_{B_\Omega(0,R)} u(y,\tau)\varphi(y,t-\tau)\,dy\\
 & \ge \int_{B_\Omega(0,R)} y_Nu(y,\tau)\phi(y)\,dy-
\sup_{x\in \Omega\cap B(0,R)}\left|x_N^{-1}\varphi(x,t-\tau)-\phi(x)\right|
\sup_{0<\tau<t}\,
\int_{B_\Omega(0,R)} y_Nu(y,\tau)\,dy\\
 & \ge \int_\Omega y_Nu(y,\tau)\phi(y)\,dy-C
\sup_{x\in \Omega\cap B(0,R)}\left|x_N^{-1}\varphi(x,t-\tau)-\phi(x)\right|
\end{split}
\end{equation*}
for all $t$, $\tau\in\Sigma_T$ with $t>\tau$. 
This together with \eqref{eq:4.13} and \eqref{eq:4.20} implies that 
\begin{align*}
\int_\Omega u(x,t)x_N\phi(x)\,dx
 & \ge\int_{\overline{\Omega}} \phi(y)\,d\nu-C\sup_{x\in \Omega\cap B(0,R)}\left|x_N^{-1}\varphi(x,t)-\phi(x)\right|\\
 & +\int_0^t\int_\Omega \int_\Omega G(x,y,t-s)u(y,s)^px_N\phi(x)\,dx\,dy\,ds.
\end{align*}
Then, by \eqref{eq:4.1} and \eqref{eq:4.20}  we see that 
$$
\underset{t\to +0}{\mbox{{\rm ess lim}}}\int_0^t\int_\Omega \int_\Omega G(x,y,t-s)u(y,s)^px_N\phi(x)\,dx\,dy\,ds=0.
$$
Since 
\begin{equation*}
\begin{split}
\int_\Omega u(x,t)x_N\phi(x)\,dx & =\int_\Omega \int_{\overline{\Omega}} K(x,y,t)x_N\phi(x)\,d\mu(y)\,dx\\
 & +\int_0^t\int_\Omega\int_\Omega G(x,y,t-s)u(y,s)^p x_N\phi(x)\,dx\,dy\,ds
\end{split}
\end{equation*}
for a.a.~$t\in(0,T)$ (see Definition~\ref{Definition:1.1}~(ii)), we have
\begin{equation}
\label{eq:4.21}
\int_{\overline{\Omega}} \phi(x)\,d\nu(x)
=\underset{t\to +0}{\mbox{{\rm ess lim}}}\int_\Omega \int_{\overline{\Omega}} K(x,y,t)x_N\phi(x)\,d\mu(y)\,dx.
\end{equation}
On the other hand, by \eqref{eq:1.5} and \eqref{eq:1.6} we see that
\begin{align*}
 & (\partial_{y_N}\varphi)(y',0,t)\\
 & =\frac{\partial}{\partial y_N}
 \int_\Omega \Gamma_{N-1}(y'-x',t)[\Gamma_1(y_N-x_N,t)-\Gamma_1(y_N+x_N)]x_N\phi(x)\,dx\biggr|_{y_N=0}\\
 & =\int_\Omega \frac{x_N}{t}\Gamma_N(x'-y',x_N,t)x_N\phi(x)\,dx=\int_\Omega K(x,y',0,t)x_N\phi(x)\,dx
\end{align*}
for all $(y',0)\in\partial\Omega$. 
Then we observe that
\begin{align*}
 & \int_\Omega \int_{\overline{\Omega}} K(x,y,t)x_N\phi(x)\,d\mu(y)\,dx\\
 & =\int_\Omega y_N^{-1}\left(\int_\Omega G(x,y,t)\varphi_0(x)\,dx\right)\,d\mu(y)
 +\int_{\partial\Omega}\left(\int_\Omega K(x,y',0,t)x_N\phi(x)\,dx\right)\,d\mu(y)\\
 & =\int_\Omega y_N^{-1}\varphi(y,t)\,d\mu(y)+\int_{\partial\Omega}(\partial_{y_N}\varphi)(y',0,t)\,d\mu(y).
\end{align*}
Therefore, by Lemmas~\ref{Lemma:4.1} and \ref{Lemma:4.3} we have 
$$
\lim_{t\to +0}\int_\Omega \int_{\overline{\Omega}} K(x,y,t)x_N\phi(x)\,d\mu(y)\,dx
=\int_{\overline{\Omega}}\phi(y)\,d\mu(y),
$$
which together with \eqref{eq:4.21} implies that
$$
\int_{\overline{\Omega}}\phi(x)\,d\nu(x)=\int_{\overline{\Omega}} \phi(y)\,d\mu(y).
$$
Since $\phi$ is arbitrary, we deduce that $\mu=\nu$ in ${\mathcal M}$. Thus assertion~(ii) follows. 
$\Box$ 
\vspace{5pt}
\newline
{\bf Proof of Theorem~\ref{Theorem:1.1}~(i).}
Let $u\in {\mathcal L}(Q_T)$  solve problem~(E) in $Q_T$, where $T>0$. 
By Lemma~\ref{Lemma:4.2} we find a unique $\nu\in{\mathcal M}$ satisfying \eqref{eq:1.10}. 
It remains to prove that $u$ is a solution to problem~(P) with $\mu=\nu$ in $Q_T$. 

Let $(x,t)\in Q_T$ be such that \eqref{eq:1.7} holds for a.a.~$\tau\in(0,T)$ with $\tau<t$. 
Similarly to the proof of Lemma~\ref{Lemma:4.3}, 
for any $n=1,2,\dots$,  
by the Besicovitch covering lemma we find an integer $m=m(N)$ and 
$\{x_{k,i}\}_{k=1,\dots,m,\,i\in{\bf N}}\subset \overline{\Omega}\setminus B(0,n\sqrt{t})$ such that 
\begin{equation}
\label{eq:4.22}
B_{k,i}\cap B_{k,j}=\emptyset\quad\mbox{if $i\not=j$}
\qquad\mbox{and}\qquad
\overline{\Omega}\setminus B(0,n\sqrt{t})\subset\bigcup_{k=1}^m\bigcup_{i=1}^\infty B_{k,i},
\end{equation}
where $B_{k,i}:=\overline{B(x_{k,i},\sqrt{t})}$. 
For any $y=(y',y_N)$, $z\in B_{k,i}$, and $\tau\in(0,t/2)$, 
since 
$$
|x-y|^2\ge (|x-z|-|z-y|)^2=|x-z|^2+|z-y|^2-2|x-z||z-y|
\ge\frac{1}{2}|x-z|^2-Ct,
$$
by Lemma~\ref{Lemma:2.2} we have 
\begin{equation}
\label{eq:4.23}
\begin{split}
K(x,y,t-\tau) & \le \frac{C}{y_N+\sqrt{t-\tau}}\Gamma_N(x-y,2(t-\tau))\\
 & \le Ct^{-\frac{N+1}{2}}\exp\left(-\frac{|x-y|^2}{Ct}\right)
\le Ct^{-\frac{N+1}{2}}\exp\left(-\frac{|x-z|^2}{Ct}\right).
\end{split}
\end{equation}
We observe from \eqref{eq:4.16}, \eqref{eq:4.22}, and \eqref{eq:4.23} that
\begin{equation*}
\begin{split}
 & \int_{\overline{\Omega}\setminus B(0,n\sqrt{t})}K(x,y,t-\tau)\,d\nu(y)
 +\int_{\Omega\setminus B(0,n\sqrt{t})}G(x,y,t-\tau)u(y,\tau)\,dy\\
 & \le\sum^m_{k=1}\sum^\infty_{i=1}\int_{B_{k,i}\cap\,\overline{\Omega}}
 K(x,y,t-\tau)(d\nu(y)+y_Nu(y,\tau)\,dy)\\
 & \le\sum^m_{k=1}\sum^\infty_{i=1}\sup_{y\in B_{k,i}}K(x,y,t-\tau)\left(\nu(B_{k,i}\cap\,\overline{\Omega})+\int_{B_{k,i}\cap\,\overline{\Omega}}y_Nu(y,\tau)\,dy\right)\\
 & \le C\sum^m_{k=1}\sum^\infty_{i=1}t^{-\frac{N+1}{2}}\inf_{z\in B_{k,i}}\exp\left(-\frac{|x-z|^2}{Ct}\right)(|x_{k,i}|+\sqrt{T})\\
 & \le Ct^{-N-\frac{1}{2}}\sum^m_{k=1}\sum^\infty_{i=1}\int_{B_{k,i}}\exp\left(-\frac{|x-z|^2}{Ct}\right)(|z|+\sqrt{T})\,dz\\
  & \le Ct^{-N-\frac{1}{2}}\int_{{\mathbb R}^N\setminus B(0,(n-1)\sqrt{t})}\exp\left(-\frac{|x-z|^2}{Ct}\right)(|z|+\sqrt{T})\,dz
\end{split}
\end{equation*}
for a.a.~$\tau\in(0,t/2)$. 
Then we obtain 
\begin{equation}
\label{eq:4.24}
\begin{split}
 & \underset{\tau\in(0,t/2)}{{\rm ess\,\,sup}}
\int_{\overline{\Omega}\setminus B(0,n\sqrt{t})}K(x,y,t-\tau)\,d\nu(y)\\
 & \qquad\quad+\underset{\tau\in(0,t/2)}{{\rm ess\,\,sup}}\,\int_{\Omega\setminus B(0,n\sqrt{t})}G(x,y,t-\tau)u(y,\tau)\,dy\to 0
 \quad\mbox{as}\quad n\to\infty.
\end{split}
\end{equation}

Let $\phi_n\in C_0^\infty({\mathbb R}^N)$ be such that 
$$
0\le\phi_n\le 1\quad\mbox{in}\quad{\mathbb R}^N,
\qquad
\phi_n=1\quad\mbox{on}\quad B(0,n\sqrt{t}),
\qquad
\phi_n=0\quad\mbox{outside}\quad B(0,2n\sqrt{t}). 
$$
Then we have 
\begin{equation}
\label{eq:4.25}
\begin{split}
 & \left|\int_\Omega G(x,y,t-\tau)u(y,\tau)\,dy-\int_{\overline{\Omega}} K(x,y,t)\,d\nu(y)\right|\\
 & \le\left|\int_\Omega G(x,y,t)u(y,\tau)\phi_n(y)\,dy - \int_{\overline{\Omega}} K(x,y,t)\phi_n(y)\,d\nu(y)\right|\\
 & \qquad
 +\left|\int_\Omega[G(x,y,t-\tau)-G(x,y,t)]u(y,\tau)\phi_n(y)\,dy\right|\\
 & \qquad\qquad
 +\int_{\Omega\setminus B(0,n\sqrt{t})}G(y,t-\tau)u(y,\tau)\,dy
 +\int_{\overline{\Omega}\setminus B(0,n\sqrt{t})} K(x,y,t)\,d\nu(y)
\end{split}
\end{equation}
for $n=1,2,\dots$ and $\tau\in(0,t/2)$. 
By \eqref{eq:1.9} we see that 
\begin{equation}
\label{eq:4.26}
\underset{t\to +0}{\mbox{{\rm ess lim}}}\,
\left[\int_\Omega G(x,y,t)u(y,\tau)\phi_n(y)\,dy - \int_{\overline{\Omega}} K(x,y,t)\phi_n(y)\,d\nu(y)\right]=0.
\end{equation}
Furthermore, by \eqref{eq:4.10} we have 
\begin{equation}
\label{eq:4.27}
\begin{split}
 & \lim_{\tau\to+0}\,\left|\int_D [G(x,y,t-\tau)-G(x,y,t)]u(y,\tau)\phi_n(y)\,dy\right|\\
 & 
\le \sup_{y\in B(0,2n\sqrt{t}),s\in(t/2,t)}\,|\partial_t K(x,y,s)|
\,\underset{\tau\to +0}{\mbox{{\rm ess limsup}}}
\biggr[\tau\int_{B_\Omega(0,2n\sqrt{t})}y_Nu(y,\tau)\,dy\biggr]=0.
\end{split}
\end{equation}
By \eqref{eq:4.25}, \eqref{eq:4.26}, and \eqref{eq:4.27} we see that 
\begin{equation*}
\begin{split}
 & \underset{\tau\to +0}{\mbox{{\rm ess limsup}}}\,
 \left|\int_\Omega G(x,y,t-\tau)u(y,\tau)\,dy - \int_{\overline{\Omega}} K(x,y,t)\,d\nu(y)\right|\\
 & \le
\underset{\tau\in(0,t/2)}{{\rm ess\,\,sup}}
\,\int_{\Omega\setminus B(0,n\sqrt{t})}G(x,y,t-\tau)u(y,\tau)\,dy
 +\int_{\overline{\Omega}\setminus B(0,n\sqrt{t})} K(x,y,t)\,d\nu(y)
\end{split}
\end{equation*}
for $n=1,2,\dots$. 
This together with \eqref{eq:4.24} implies that 
$$
\underset{\tau\to +0}{\mbox{{\rm ess lim}}}\,
\left|\int_\Omega G(x,y,t-\tau)u(y,\tau)\,dy - \int_{\overline{\Omega}} K(x,y,t)\,d\nu(y)\right|=0.
$$
Therefore we observe from Definition~\ref{Definition:1.1}~(i) that 
$u$ is a solution to problem~(P) with $\mu=\nu$ in~$Q_T$. 
Thus Theorem~\ref{Theorem:1.1}~(ii) follows. The proof of Theorem~\ref{Theorem:1.1} is complete. 
$\Box$
\medskip

Finally, we prove Corollaries~\ref{Corollary:1.1} and \ref{Corollary:1.2}.
\vspace{5pt} 
\newline
{\bf Proof of Corollary~\ref{Corollary:1.1}.}
Corollary~\ref{Corollary:1.1} follows from Theorem~\ref{Theorem:1.1} and Theorem~\ref{Theorem:1.2}~(iii). 
$\Box$
\vspace{5pt} 
\newline
{\bf Proof of Corollary~\ref{Corollary:1.2}.}
By Lemma~\ref{Lemma:2.4} we see that, 
for a.a.~$\tau\in(0,T)$, 
$u_\tau(x,t):=u(x,t+\tau)$ is a solution to problem~\eqref{eq:P} in $Q_{T-\tau}$ 
with $\mu(x)=x_Nu(x,\tau)$ on $\overline{\Omega}$. 
Then Corollary~\ref{Corollary:1.2}. follows from Theorem~\ref{Theorem:1.2}. 
$\Box$
%
\section{Sufficient conditions}
In this section we study sufficient conditions for the solvability of problem~\eqref{eq:P}. 
We denote by ${\mathcal L}$ (resp.\,\,${\mathcal L}'$) 
the set of nonnegative measurable functions in $\Omega$ (resp.\,\,$\partial\Omega$). 
For any $f\in {\mathcal L}$ and $h\in {\mathcal L'}$, 
we set 
$$
[G(t)f](x):=\int_\Omega G(x,y,t)f(y)\,dy, \quad 
[\Gamma_{N-1}(t)h](x'):=\int_{\partial \Omega} \Gamma_{N-1}(x',y',t)h(y')\,dy', 
$$
for $x\in\Omega$, $x'\in{\mathbb R}^{N-1}$, and $t>0$. 
\subsection{Sufficient conditions in the case of $\mu\in{\mathcal M}$}
We start by showing the following theorem.
\begin{theorem}
\label{Theorem:5.1}
Let $N\ge 1$ and $p>1$. 
Then there exists $\gamma=\gamma(N,p)>0$ such that, 
if $\mu\in{\mathcal M}$ satisfies
\begin{equation}
\label{eq:5.1}
\int_0^T s^{-\frac{N(p-1)}{2}}
\left(\sup_{z\in{\overline{\Omega}}}\int_{B_\Omega(z,\sqrt{s})}\frac{d\mu(y)}{y_N+\sqrt{s}}\right)^{p-1}\,ds\le \gamma
\end{equation}
for some $T>0$, then problem~\eqref{eq:P} possesses a solution in $Q_T$. 
\end{theorem}
{\bf Proof.}
Assume \eqref{eq:5.1}. Let $T>0$ and 
$$
w(x,t)=2\int_{\overline{\Omega}}K(x,y,t)\,d\mu(y). 
$$
It follows from Lemma~\ref{Lemma:2.2} that
$$
\|w(t)\|_{L^\infty(\Omega)}\le Ct^{-\frac{N}{2}}\sup_{z\in{\overline{\Omega}}}\int_{B_\Omega(z,\sqrt{t})}\frac{d\mu(y)}{y_N+\sqrt{t}}
\quad\mbox{for all $t>0$}.
$$
Then, by Lemma~\ref{Lemma:2.3} and \eqref{eq:5.1} we have
\begin{align*}
 & \int_{\overline{\Omega}}K(x,y,t)\,d\mu(y)+\int_0^t\int_\Omega G(x,y,t-s)w(y,s)^p\,dy\,ds\\
 & \le\frac{1}{2}w(x,t)\\
 & +C\int_0^t
 s^{-\frac{N(p-1)}{2}}\left(\sup_{z\in{\overline{\Omega}}}\int_{B_\Omega(z,\sqrt{s})}\frac{d\mu(y)}{y_N+\sqrt{s}}\right)^{p-1}
 \left(\int_{\overline{\Omega}}\int_\Omega G(x,y,t-s)K(y,z,s)\,dy\,d\mu(z)\right)\,ds\\
  & \le\frac{1}{2}w(x,t)+C\gamma w(x,t)
\end{align*}
for a.a.~$(x,t)\in Q_T$. 
If $\gamma>0$ is small enough, then we see that $w$ is a supersolution to problem~\eqref{eq:P} in $Q_T$. 
Thus Theorem~\ref{Theorem:5.1} follows from Lemma~\ref{Lemma:2.5}.
$\Box$\vspace{5pt}\newline
As corollaries of Theorems~\ref{Theorem:1.2} and \ref{Theorem:5.1}, we have 
the following results.
\begin{corollary}
\label{Corollary:5.1}
Let $N\ge 1$ and $\mu\in{\mathcal M}$. 
\begin{itemize}
  \item[{\rm (i)}] 
  Let $1<p<p_{N+1}$. 
  Then problem~\eqref{eq:P} possesses a local-in-time solution if and only if 
  \begin{equation}
  \label{eq:5.2}
  \sup_{z=(z',z_N)\in\overline{\Omega}}\frac{\mu(B_\Omega(z,1))}{1+z_N}<\infty. 
  \end{equation}
  \item[{\rm (ii)}] 
  If $\mbox{supp}\,\mu\subset\Omega_L$ for some $L>0$, then assertion~{\rm (i)} holds for $1<p<p_N$. 
\end{itemize}
\end{corollary}
{\bf Proof.}
Let $p>1$. 
Assume that problem~\eqref{eq:P} possesses a solution in $Q_T$ for some $T>0$. 
It follows from Theorem~\ref{Theorem:1.2} that 
$$
\mu(B_\Omega(z,\sqrt{T}))=\lim_{\sigma\to\sqrt{T}}\mu(B_\Omega(z,\sigma))
\le CT^{N-\frac{1}{p-1}}(z_N+\sqrt{T})
$$
for all $z\in\overline{\Omega}$. This together with Lemma~\ref{Lemma:2.1} implies \eqref{eq:5.2}. 

Conversely, assume that \eqref{eq:5.2} holds. 
Consider the case of $1<p<p_{N+1}$. 
Since 
$$
\int_{B_\Omega(z,\sqrt{s})}\frac{d\mu(y)}{y_N+\sqrt{s}}\le
\left\{
\begin{array}{ll}
s^{-\frac{1}{2}}\mu(B_\Omega(z,\sqrt{s})) & \quad\mbox{if}\quad z_N\le 2,\vspace{5pt}\\
2z_N^{-1}\mu(B_\Omega(z,\sqrt{s}))& \quad\mbox{if}\quad z_N>2,
\end{array}
\right.
$$
for all $z\in\overline{\Omega}$ and $s\in(0,1)$, by \eqref{eq:5.2} we have
$$
\int_0^T s^{-\frac{N(p-1)}{2}}
\left(\sup_{z\in{\overline{\Omega}}}\int_{B_\Omega(z,\sqrt{s})}\frac{d\mu(y)}{y_N+\sqrt{s}}\right)^{p-1}\,ds
\le C\int_0^T s^{-\frac{(N+1)(p-1)}{2}}\,ds
\le CT^{1-\frac{(N+1)(p-1)}{2}}
$$
for all $T\in(0,1)$. Then Theorem~\ref{Theorem:5.1} implies that 
problem~\eqref{eq:P} possesses a local-in-time solution. 
Thus assertion~(i) follows. 

Next, we assume that $\mbox{supp}\,\mu\subset\Omega_L$ for some $L>0$, and consider the case of $1<p<p_N$. 
Assume \eqref{eq:5.2}. Then
$$
\int_{B_\Omega(z,\sqrt{s})}\frac{d\mu(y)}{y_N+\sqrt{s}}\le
C(1+z_N)^{-1}\mu(B_\Omega(z,\sqrt{s}))\le C
$$
for all $z\in\overline{\Omega}$ and small enough $s>0$. 
This together with $1<p<p_N$ implies that  
$$
\int_0^T s^{-\frac{N(p-1)}{2}}
\left(\sup_{z\in{\overline{\Omega}}}\int_{B_\Omega(z,\sqrt{s})}\frac{d\mu(y)}{y_N+\sqrt{s}}\right)^{p-1}\,ds
\le C\int_0^T s^{-\frac{N(p-1)}{2}}\,ds
\le CT^{1-\frac{N(p-1)}{2}}
$$
for all small enough $T>0$. 
Then we observe from Theorem~\ref{Theorem:5.1} that 
problem~\eqref{eq:P} possesses a local-in-time solution, 
and assertion~(ii) follows from assertion~(i). 
Thus Corollary~\ref{Corollary:5.1} follows.~$\Box$
\begin{corollary}
\label{Corollary:5.2}
Let $N\ge 1$. Problem~{\rm (E)} possesses a nontrivial global-in-time solution if and only if $p>p_{N+1}$. 
\end{corollary}
{\bf Proof.}
Assume that there exists a nontrivial global-in-time solution $u$ to problem~\eqref{eq:P}. 
Then, by Lemma~\ref{Lemma:2.4} we find $\tau\in(0,\infty)$ such that 
$u_\tau(x,t):=u(x,t+\tau)$ is a global-in-time solution 
to problem~\eqref{eq:P} with $\mu=x_Nu(x,\tau)$ on $\overline{\Omega}$ 
and $u(\cdot,\tau)\not=0$ in ${\mathcal L}$.  
Then, by Theorem~\ref{Theorem:1.2} we have 
$$
0\le\int_{B_\Omega(0,\sigma)}y_Nu(y,\tau)\,dy\le C\sigma^{-\frac{2}{p-1}}\int_{B_\Omega(0,\sigma)}y_N\,dy
\le C\sigma^{N+1-\frac{2}{p-1}}
\quad\mbox{for all $\sigma>0$}. 
$$
If $1<p<p_{N+1}$, then
$$
0\le\int_\Omega y_Nu(y,\tau)\,dy\le C\lim_{\sigma\to\infty}\sigma^{N+1-\frac{2}{p-1}}=0.
$$
This is a contradiction. 
If $p=p_{N+1}$, then, by Theorem~\ref{Theorem:1.2}~(ii) with $\sigma\ge 1$ and $T=\sigma^4$ 
we have 
$$
0\le\int_\Omega y_Nu(y,\tau)\,dy\le C\lim_{\sigma\to\infty}\left[\log\left(e+\frac{\sqrt{\sigma^4}}{\sigma}\right)\right]^{-\frac{N+1}{2}}=0.
$$
This is also a contradiction. Thus problem~{\rm (E)} possesses no nontrivial global-in-time solutions if $1<p\le p_{N+1}$. 
On the other hand, by Theorem~\ref{Theorem:5.1}
we easily find a global-in-time positive solution to problem~\eqref{eq:P} if $p>p_{N+1}$. 
Then we complete the proof of Corollary~\ref{Corollary:5.2}.
$\Box$
\begin{remark}
\label{Remark:5.1} 
It has already been proved in \cites{LM, Mei} that the case $p=p_{N+1}$ is the threshold 
for the existence of nontrivial global-in-time solution to problem~\eqref{eq:E}. 
{\rm ({\it See also} \cites{DL, Levine} {\it for related results}.)} 
On the other hand, the proof of Corollary~{\rm\ref{Corollary:5.2}} asserts that 
our necessary conditions and sufficient conditions are useful even for the study of the existence of global-in-time solutions. 
\end{remark}
\begin{corollary}
\label{Corollary:5.3}
Let $N\ge 1$ and $\kappa>0$. 
Consider problem~\eqref{eq:P} with $\mu=\kappa\delta_N$ on $\overline{\Omega}$, that is, 
$$
u(x,0)=-\kappa\,\partial_{x_N}\delta_N\quad\mbox{on}\quad\overline{\Omega}.
$$
\begin{itemize}
  \item[{\rm (i)}] 
  If $p\ge p_{N+1}$, then problem~\eqref{eq:P} possesses no local-in-time solutions. 
  \item[{\rm (ii)}] 
  If $1<p<p_{N+1}$, then 
  problem~\eqref{eq:P} possesses a local-in-time solution. 
\end{itemize}
\end{corollary}
{\bf Proof.}
Assertion~(i) follows from Theorem~\ref{Theorem:1.2}. (See \eqref{eq:1.10} and \eqref{eq:1.12}.)
Assertion~(ii) is proved by Corollary~\ref{Corollary:5.1}. Thus Corollary~\ref{Corollary:5.3} follows. 
$\Box$
\subsection{More delicate sufficient conditions}
In this subsection 
we modify the arguments in \cites{FHIL, HI18, RS} 
to obtain the following theorem on sufficient conditions for the solvability of problem~\eqref{eq:P}. 
We denote by $\|\cdot\|_\infty$ the usual norm of $L^\infty(\Omega)$. 
\begin{theorem}
\label{Theorem:5.2}
Let $f\in{\mathcal L}$ and let $h\in{\mathcal L}'$ if $1<p<2$ and $h=0$ if $p\ge 2$. 
Consider problem~\eqref{eq:P} with 
\begin{equation}
\label{eq:5.3}
\mu=x_Nf(x)+h(x')\otimes\delta_1(x_N)\in{\mathcal M}.
\end{equation}
Let $\Phi$ be a strictly increasing, nonnegative, and convex function on $[0,\infty)$. Set
\[
v(x,t):=2\Phi^{-1}\left([G(t)\Phi(f)](x)\right),\quad
w(x,t):=2\frac{x_N}{t} \Gamma_1(x_N,t) 
\Phi^{-1}\left([\Gamma_{N-1}(t)\Phi(h)](x)\right),
\]
for $(x,t)\in Q_\infty$. Define
\[
	A(\tau):=\frac{\Phi^{-1}(\tau)^p}{\tau},\quad
	B(\tau):=\frac{\tau}{\Phi^{-1}(\tau)},\quad\mbox{for $\tau>0$}. 
\]
If 
\begin{equation}
\label{eq:5.4}
\begin{split}
 & \sup_{t\in(0,T)}
\left(\|B( G(t) \Phi(f) )\|_\infty
\int_0^t \|A( G(s) \Phi(f) )\|_\infty\,ds\right)
\le 2^{-2p+1},\\
 & \sup_{t\in(0,T)}
\left(\|B( \Gamma_{N-1}(t)\Phi(h) )\|_\infty
\int_0^t s^{-(p-1)} \| A( \Gamma_{N-1}(s) \Phi(h) )\|_\infty \,ds \right)
\le 2^{-2p+1}(2e\pi)^\frac{p-1}{2},
\end{split}
\end{equation}
for some $T>0$, then problem~\eqref{eq:P} possesses a solution $u$ in $Q_T$ such that 
$$
0\le u(x,t)\le v(x,t)+w(x,t)\quad\mbox{for a.a.~$(x,t)\in Q_T$}.
$$
\end{theorem}
{\bf Proof.}
Let $\mu$ be as in \eqref{eq:5.3}. 
We show that $v+w$ is a supersolution to problem~\eqref{eq:P} in $Q_T$. 
By \eqref{eq:1.5}, \eqref{eq:1.6}, and \eqref{eq:5.3} we have
\begin{equation}
\label{eq:5.5}
\begin{split}
[K(t)\mu](x):= &\,\int_{\overline{\Omega}}K(x,y,t)\,d\mu(y)\\
 = &\,\int_\Omega G(x,y,t)f(y)\,dy+\frac{x_N}{t}\int_{{\mathbb R}^{N-1}}\Gamma_{N-1}(x'-y',t)h(y')\,dy'
\end{split}
\end{equation}
for all $(x,t)\in Q_\infty$.
It follows from Jensen's inequality with the convexity of $\Phi$ 
and~\eqref{eq:5.5} that 
$$
[K(t)\mu](x)
 \leq \Phi^{-1}([G(t)\Phi(f)](x))
+\frac{x_N}{t} \Gamma_1(x_N,t) \Phi^{-1}([\Gamma_{N-1}(t)\Phi(h)](x'))  
= \frac{v(x,t)+w(x,t)}{2}
$$
for all $(x,t)\in Q_\infty$.
Since $(a+b)^p\leq 2^{p-1}(a^p+b^p)$ for $a,b>0$, we have 
\[
\begin{aligned}
	& K(t)\mu+\int_0^t G(t-s)(v(s)+w(s))^p\,ds \\
	&\leq 
	\frac{v(t)+w(t)}{2}+2^{p-1}\int_0^t G(t-s)v(s)^p\,ds
	+ 2^{p-1} \int_0^t G(t-s)w(s)^p\,ds. 
\end{aligned}
\]
By the semigroup property of $G$ and \eqref{eq:5.4} we see that 
\[
\begin{aligned}
	\int_0^t G(t-s)v(s)^p\,ds &
	\leq 2^p \int_0^t G(t-s) 
	\left\| \frac{[\Phi^{-1}(G(s)\Phi(f))]^p}{G(s)\Phi(f)} \right\|_{\infty}
	G(s) \Phi(f) \,ds \\
	&=
	2^p G(t)\Phi(f) \int_0^t 
	\left\| \frac{[\Phi^{-1}(G(s)\Phi(f))]^p}{G(s)\Phi(f)} \right\|_{\infty}
	\,ds \\
	&\leq 
	2^{p-1} v(t)
	\left\| \frac{G(t)\Phi(f)}{\Phi^{-1}(G(t)\Phi(f))}\right\|_\infty 
	\int_0^t 
	\left\| \frac{[\Phi^{-1}(G(s)\Phi(f))]^p}{G(s)\Phi(f)} \right\|_{\infty}
	\,ds \leq \frac{v(t)}{2^p}. 
\end{aligned}
\]
On the other hand, since
\[
	\frac{y_N}{s} \Gamma_1(y_N,s) 
	\leq 
	\frac{1}{s} (4\pi)^{-\frac{1}{2}} 
	\sup_{y_N>0} \frac{y_N}{s^{1/2}} e^{-\frac{y_N^2}{4s}} 
	\leq (2e\pi)^{-\frac{1}{2}} s^{-1}
	\quad\mbox{for all $(y_N,s)\in(0,\infty)^2$}, 
\]
by the semigroup property of $\Gamma_{N-1}$ we have
\[
\begin{aligned}
	&\int_0^t G(t-s)w(s)^p\,ds \\
	&=
	2^p \int_0^t \int_0^\infty 
	[\Gamma_1(x_N-y_N,t-s)-\Gamma_1(x_N+y_N,t-s)] 
	\left( \frac{y_N}{s} \Gamma_1(y_N,s) \right)^p dy_N \\
	&\qquad 
	\times \int_{{\mathbb R}^{N-1}} \Gamma_{N-1}(x'-y',t-s) 
	\frac{[\Phi^{-1}([\Gamma_{N-1}(s)\Phi(h)](y'))]^p} 
	{[\Gamma_{N-1}(s)\Phi(h)](y')} [\Gamma_{N-1}(s)\Phi(h)](y')
	\,dy'\,ds \\
	&\leq 
	c_p \int_0^t 
	s^{-(p-1)} 
	\int_0^\infty 
	[\Gamma_1(x_N-y_N,t-s)-\Gamma_1(x_N+y_N,t-s)] 
	\frac{y_N}{s} \Gamma_1(y_N,s)\,dy_N \\
	&\qquad 
	\times  \left\| \frac{[\Phi^{-1}([\Gamma_{N-1}(s)\Phi(h)])]^p}
	{\Gamma_{N-1}(s)\Phi(h)} \right\|_\infty 
	\int_{{\mathbb R}^{N-1}} \Gamma_{N-1}(x'-y',t-s) 
	[\Gamma_{N-1}(s)\Phi(h)](y')
	\,dy'\,ds \\
	&\leq 
	c_p \int_0^t 
	s^{-(p-1)} 
	\int_0^\infty 
	[\Gamma_1(x_N-y_N,t-s)-\Gamma_1(x_N+y_N,t-s)] 
	\frac{y_N}{s} \Gamma_1(y_N,s)\,dy_N \\
	&\qquad 
	\times  \left\| \frac{[\Phi^{-1}([\Gamma_{N-1}(s)\Phi(h)])]^p}
	{\Gamma_{N-1}(s)\Phi(h)} \right\|_\infty 
	\left\| \frac{\Gamma_{N-1}(t) \Phi(h)}
	{\Phi^{-1}(\Gamma_{N-1}(t) \Phi(h))} \right\|_\infty
	\Phi^{-1} ( [\Gamma_{N-1}(t) \Phi(h)](x') )
	\,ds, 
\end{aligned}
\]
where $c_p:=2^p(2e\pi)^{-\frac{p-1}{2}}$. 
Since $(y_N/s) \Gamma_1(y_N,s)=-2\partial_{y_N} \Gamma_1(y_N,s)$, 
integration by parts and the semigroup property of $\Gamma_1$ show that 
\[
\begin{aligned}
	&\int_0^\infty 
	[\Gamma_1(x_N-y_N,t-s)-\Gamma_1(x_N+y_N,t-s)] 
	\frac{y_N}{s} \Gamma_1(y_N,s)\,dy_N \\
	&=
	2\int_0^\infty \partial_{y_N} [ \Gamma_1(x_N-y_N,t-s)-\Gamma_1(x_N+y_N,t-s) ] 
	\Gamma_1(y_N,s)\,dy_N \\
	&=
	-2 \partial_{x_N} \int_0^\infty [ \Gamma_1(x_N-y_N,t-s)+\Gamma_1(x_N+y_N,t-s) ]
	\Gamma_1(y_N,s)\,dy_N \\
	&= -2 \partial_{x_N} \Gamma_1(x_N,t) 
	= \frac{x_N}{t} \Gamma_1(x_N,t). 
\end{aligned}
\]
Then, by \eqref{eq:5.4} 
we see that 
\[
\begin{aligned}
	&\int_0^t G(t-s)w(s)^p\,ds \\
	&\leq 
	2^p (2e\pi)^{-\frac{p-1}{2}} 
	\frac{x_N}{t} \Gamma_1(x_N,t) 
	\Phi^{-1} ( [\Gamma_{N-1}(t) \Phi(h)](x') )  \\
	&\qquad \quad \times 
	\left\| \frac{\Gamma_{N-1}(t) \Phi(h)}
	{\Phi^{-1}(\Gamma_{N-1}(t) \Phi(h))} \right\|_\infty
	\int_0^t s^{-(p-1)} \left\| \frac{[\Phi^{-1}([\Gamma_{N-1}(s)\Phi(h)])]^p}
	{\Gamma_{N-1}(s)\Phi(h)} \right\|_\infty\,ds \\
	&\leq 
	2^{p-1} (2e\pi)^{-\frac{p-1}{2}} w(t)
	\| B( \Gamma_{N-1}(t) \Phi(h) ) \|_\infty 
	\int_0^t  s^{-(p-1)}\| A( \Gamma_{N-1}(s) \Phi(h) ) \|_\infty\,ds
	\leq 
	\frac{w(t)}{2^p}.
\end{aligned}
\]
The above computations show that 
\[
 K(t)\mu+\int_0^t G(t-s)(v(s)+w(s))^p\,ds
 \le v(t)+w(t)\quad\mbox{for all $t\in(0,T)$}.
\]
This means that $v+w$ is a supersolution to problem~\eqref{eq:P} in $Q_T$. 
Then Lemma~\ref{Lemma:2.5} implies that problem~\eqref{eq:P} possesses a solution in 
$Q_T$. 
Thus Theorem~\ref{Theorem:5.2} follows.
$\Box$\vspace{5pt}

Next, as an application of Theorem~\ref{Theorem:5.2},
we obtain sufficient conditions for the solvability of problem~\eqref{eq:P}. 
\begin{theorem}
\label{Theorem:5.3}
Let $f\in{\mathcal L}$ and let $h\in{\mathcal L}'$ if $1<p<2$ and $h=0$ if $p\ge 2$. 
For any $\alpha>1$, there exists $\gamma=\gamma(N,p,\alpha)>0$ with the following property: 
if there exists $T>0$ such that
\begin{equation}
\label{eq:5.6}
\begin{split}
 & \sup_{x\in\overline{\Omega}}\int_{B_\Omega(x,\sigma)}
 \frac{y_N}{y_N+\sigma}f(y)^\alpha\,dy
\le\gamma \sigma^{N-\frac{2\alpha}{p-1}},\\
 & \sup_{x'\in{\bf R}^{N-1}}\int_{B'(x',\sigma)}h(y')^\alpha\,dy'\le\gamma\sigma^{N-1+2\alpha\frac{p-2}{p-1}},
\quad\mbox{for all $\sigma\in(0,\sqrt{T})$},
\end{split}
\end{equation}
then problem~\eqref{eq:P} with \eqref{eq:5.3} 
possesses a solution~$u$ in $Q_T$, with $u$ satisfying 
\begin{equation}
\label{eq:5.7}
0\le u(x,t)\le 2[G(t)f^\alpha](x)^{\frac{1}{\alpha}}
+2\frac{x_N}{t}\Gamma_1(x_N,t) [\Gamma_{N-1}(t)h^\alpha](x')^{\frac{1}{\alpha}}
\quad\mbox{for a.a.~$(x,t)\in Q_T$}.
\end{equation}
\end{theorem}
{\bf Proof.}
Assume \eqref{eq:5.6}. 
We can assume, without loss of generality, that $\alpha\in(1,p)$. 
Indeed, if~$\alpha\ge p$, then, for any $1<\alpha'<p$, we apply H\"older's inequality to obtain
\begin{align*}
&\sup_{x\in\overline{\Omega}}\int_{B_\Omega(x,\sigma)}
\frac{y_N}{y_N+\sigma}f(y)^{\alpha'}\,dy \\
&\le\sup_{x\in{\overline{\Omega}}}\left[\,\int_{B_\Omega(x,\sigma)}
\frac{y_N}{y_N+\sigma}\,dy\,\right]^{1-\frac{\alpha'}{\alpha}}\left[\,\int_{B_\Omega(x,\sigma)}
\frac{y_N}{y_N+\sigma}f(y)^\alpha\,dy\,\right]^{\frac{\alpha'}{\alpha}}
\le C\gamma^\frac{\alpha'}{\alpha} \sigma^{N-\frac{2\alpha'}{p-1}},\\
 & \sup_{x'\in{\bf R}^{N-1}}\int_{B'(x',\sigma)}h(y')^{\alpha'} \,dy'
\le \sup_{x'\in{\bf R}^{N-1}}\left[\,\int_{B'(x',\sigma)} \,dy' \,\right]^{1-\frac{\alpha'}{\alpha}}\left[\,\int_{B'(x',\sigma)}h(y')^{\alpha} \,dy' \,\right]^\frac{\alpha'}{\alpha} \\
&\le C\gamma^\frac{\alpha'}{\alpha}\sigma^{N-1-2\alpha'\frac{2-p}{p-1}}
\end{align*}
for all $\sigma\in(0,\sqrt{T})$. 
Thus \eqref{eq:5.6} holds with $\alpha$ replaced by $\alpha'$. 
Furthermore, if \eqref{eq:5.7} holds for some $\alpha'\in(1,\alpha)$, then, 
since 
$$
[G(t)f^{\alpha'}](x)^{\frac{1}{\alpha'}}\le [G(t)f^{\alpha}](x)^{\frac{1}{\alpha}},
\qquad
[\Gamma_{N-1}(t)h^{\alpha'}](x)^{\frac{1}{\alpha'}}\le[\Gamma_{N-1}(t)h^{\alpha}](x)^{\frac{1}{\alpha}},
$$
for $x\in\Omega$, $x'\in{\mathbb R}^{N-1}$, and $t>0$,  the desired inequality \eqref{eq:5.7} holds. 

We apply Theorem~\ref{Theorem:5.2} to prove Theorem~\ref{Theorem:5.3}. 
Let $A$ and $B$ be as in Theorem~\ref{Theorem:5.2} with $\Phi(\tau)=\tau^\alpha$. 
Then 
$A(\tau)=\tau^{\frac{p}{\alpha}-1}$ and $B(\tau)=\tau^{1-\frac{1}{\alpha}}$.
Set 
\[
v(x,t):=2[G(t)f^\alpha](x)^{\frac{1}{\alpha}},
\quad
w(x,t):=2\frac{x_N}{t}\Gamma_1(x_N,t) [\Gamma_{N-1}(t)h^\alpha](x')^{\frac{1}{\alpha}},
\]
for all $(x,t)\in Q_T$. 
It follows from Lemma~\ref{Lemma:2.2} and~\eqref{eq:5.6} that 
\begin{equation*}
\begin{split}
 & [G(t)f^\alpha](x)=\int_\Omega K(x,y,t)y_Nf(y)^\alpha\,dy
\le Ct^{-\frac{N}{2}}\sup_{x\in\overline{\Omega}}\int_{B_\Omega(x,\sqrt{t})}
\frac{y_N}{y_N+\sqrt{t}}f(y)^\alpha\,dy
\le C\gamma t^{-\frac{\alpha}{p-1}},\\
 & [\Gamma_{N-1}(t)h^\alpha](x)
\leq Ct^{-\frac{N-1}{2}} \sup_{x'\in{\mathbb R}^{N-1}} 
\int_{B'(x',\sqrt{t})} h(y')^\alpha dy'
\le C\gamma t^{\alpha\left(1-\frac{1}{p-1}\right)},
\end{split}
\end{equation*}
for all $t\in(0,\sqrt{T})$. 
Then, thanks to $\alpha\in(1,p)$, we have
\[
	\|B( G(t) f^\alpha )\|_\infty
	\int_0^t \|A( G(s) f^\alpha )\|_\infty\,ds 
	= \| G(t)f^\alpha\|_\infty^{1-\frac{1}{\alpha}}
	\int_0^t \| G(s) f^\alpha \|_\infty^{\frac{p}{\alpha}-1} \,ds 
	\leq C\gamma^\frac{p-1}{\alpha}
\]
for all $t\in(0,T)$. 
In the case of $1<p<2$, we obtain 
\[
\begin{aligned}
	&\|B( \Gamma_{N-1}(t)h^\alpha )\|_\infty
	\int_0^t s^{-(p-1)} \| A( \Gamma_{N-1}(s) h^\alpha )\|_\infty \,ds  \\
	&=
	\| \Gamma_{N-1}(t)h^\alpha\|_\infty^{1-\frac{1}{\alpha}}
	\int_0^t s^{-(p-1)} 
	\| \Gamma_{N-1}(s) h^\alpha \|_\infty^{\frac{p}{\alpha}-1} \,ds \\
	&\leq 
	C\gamma^\frac{p-1}{\alpha} t^{\left( 1-\frac{1}{p-1}\right)(\alpha-1)}
	\int_0^t s^{1-\alpha-\frac{p-\alpha}{p-1}} \,ds 
	= C\gamma^\frac{p-1}{\alpha}
	\quad\mbox{for all $t\in(0,T)$}.
\end{aligned}
\]
Then we apply Theorem~\ref{Theorem:5.2} to obtain the desired conclusion. The proof is complete. 
$\Box$
\begin{theorem}
\label{Theorem:5.4}
Let $p=p_{N+\ell}$ with $\ell\in[0,1]$. 
Let $\beta>0$ and set
$\Phi(\tau):=\tau[\log(e+\tau)]^\beta$ for $\tau\ge 0$. 
For any $T>0$, there exists $\gamma=\gamma(N,\beta,T,\ell)>0$ such that, 
if $f\in{\mathcal L}$ satisfies 
\begin{equation*}
\sup_{x\in\overline{\Omega}}\int_{B_\Omega(x,\sigma)}y_N^\ell\Phi(T^{\frac{1}{p-1}}f(y))\,dy
\le\gamma T^{\frac{N+\ell}{2}}\biggr[\log\left(e+\frac{\sqrt{T}}{\sigma}\right)\biggr]^{\beta-\frac{N+\ell}{2}}\quad
\mbox{for all $\sigma\in(0,\sqrt{T})$},
\end{equation*}
then problem~\eqref{eq:P} with $\mu=x_Nf(x)$ 
possesses a solution~$u$ in $Q_T$, with $u$ satisfying 
\[
0\le u(x,t)\le C\Phi^{-1}\left([G(t)\Phi(T^{\frac{1}{p-1}}f)](x)\right)\quad \mbox{for a.a.~$(x,t)\in Q_T$},
\]
for some $C>0$. 
\end{theorem}
{\bf Proof.}
Let $0<\epsilon<p-1$. 
We find $L\in[e,\infty)$ with the following properties: 
\begin{itemize}
\item[(a)] $\Psi(s):= s[\log(L+s)]^\beta$ is positive and convex in $(0,\infty)$;
\item[(b)] $s^p/\Psi(s)$ is increasing in $(0,\infty)$;
\item[(c)] $s^{\epsilon}[\log(L+s)]^{-\beta p}$ is increasing in $(0,\infty)$. 
\end{itemize}
Since $C^{-1}\Phi(s)\leq \Psi(s)\leq C \Phi(s)$ for $s\in(0,\infty)$, 
we see that
\begin{equation}
\label{eq:5.8}
\sup_{x\in\overline{\Omega}}\int_{B_\Omega(x,\sigma)}y_N^\ell\Psi(T^{\frac{1}{p-1}}f(y))\,dy\\
\le C\gamma T^{\frac{N+\ell}{2}}\biggr[\log\left(e+\frac{\sqrt{T}}{\sigma}\right)\biggr]^{\beta-\frac{N+\ell}{2}}
\end{equation}
for all $\sigma\in(0,\sqrt{T})$. 
Here we can assume, without loss of generality, that $\gamma\in(0,1)$. 
Set
\[
z(x,t) :=\left[G(t)\Psi(T^{\frac{1}{p-1}}f)\right](x)
=\int_\Omega K(x,y,t)y_N\Psi(T^{\frac{1}{p-1}}f(y))\,dy.
\]
By Lemma~\ref{Lemma:2.2} and \eqref{eq:5.8} we have 
\begin{equation*}
\begin{split}
\|z(t)\|_\infty & \le C t^{-\frac{N}{2}} \sup_{x\in \overline{\Omega}} \int_{B_\Omega(x,\sqrt{t} )}\frac{y_N}{y_N+\sqrt{t}}\Psi(T^{\frac{1}{p-1}}f)\,dy\\
 & \le C t^{-\frac{N+\ell}{2}} \sup_{x\in \overline{\Omega}} \int_{B_\Omega(x,\sqrt{t} )}y_N^\ell\Psi(T^{\frac{1}{p-1}}f)\,dy\\
 &\le C\gamma t_T^{-\frac{N+\ell}{2}}|\log t_T|^{\beta-\frac{N+\ell}{2}} \le C t_T^{-\frac{N+\ell}{2}}|\log t_T|^{\beta-\frac{N+\ell}{2}}
\end{split}
\end{equation*}
for all $t\in(0,T)$, where $t_T:=t/(2T)\in(0,1/2)$. 
Since 
\[
C^{-1}\tau[\log(L+\tau)]^{-\beta}\le\Psi^{-1}(\tau)
\le C\tau[\log(L+\tau)]^{-\beta}\quad\mbox{for $\tau>0$},
\]
we have
\[
\begin{aligned}
	&A(z(x,t)) = \frac{\Psi^{-1}(z(x,t))^p}{z(x,t)} 
	\leq Cz(x,t)^{p-1} [\log(L+z(x,t))]^{-\beta p}, \\
	&B(z(x,t)) = \frac{z(x,t)}{\Psi^{-1}(z(x,t))} 
	\leq C [\log(L+z(x,t))]^\beta,
\end{aligned}
\]
for $(x,t)\in Q_\infty$. 
Then we have 
\[
\begin{aligned}
	&\begin{aligned}
	0&\leq A(z(x,t))  \leq 
	C\|z(t)\|_\infty^{p-1-\epsilon} \|z(t)\|_\infty^{\epsilon} [\log(L+\|z(t)\|_\infty)]^{-\beta p}\\
	&\leq 
	C\gamma^{p-1-\epsilon} t_T^{-\frac{(N+\ell)(p-1)}{2}}
	|\log t_T|^{(\beta-\frac{N+\ell}{2})(p-1)}
	|\log t_T|^{-\beta p} 
	= C\gamma^{p-1-\epsilon} t_T^{-1}|\log t_T|^{-\beta-1}, 
	\end{aligned} \\
	&0\leq B(z(x,t)) 
	\leq C[\log(L+\|z(t)\|_\infty)]^\beta 
	\leq C|\log t_T|^\beta, 
\end{aligned}
\]
for all $(x,t)\in Q_T$, where $C$ is independent of $\gamma$.
Hence 
\[
\begin{aligned}
	&\|B(z(t))\|_\infty
	\int_0^t \|A(z(s))\|_\infty\,ds 
	\le C\gamma^{p-1-\epsilon}|\log t_T|^\beta\int_0^t s_T^{-1}|\log s_T|^{-\beta-1}\,ds \\
	&=C\gamma^{p-1-\epsilon}|\log t_T|^\beta \int_0^t \frac{2T}{s} 
	\left[ -\log \frac{s}{2T} \right]^{-\beta-1}\,ds \\
	&=C\gamma^{p-1-\epsilon}|\log t_T|^\beta 
	\left. \left( \frac{2T}{\beta} \left[ -\log \frac{s}{2T} \right]^{-\beta} 
	\right) \right|^{s=t}_{s=0}
	= C T \gamma^{p-1-\epsilon}
\end{aligned}
\]
for all $t\in(0,T)$. 
Therefore, if $\gamma>0$ is small enough, 
then we apply Theorem~\ref{Theorem:5.2} to find a solution~$u$ to problem~\eqref{eq:P} in $Q_T$ 
such that 
\[
0\le u(x,t)\le 2\Psi^{-1}(z(x,t))
\le C\Phi^{-1}([G(t)\Phi(f)](x))\quad\mbox{for a.a.~$(x,t)\in Q_\infty$}.
\]
Thus Theorem~\ref{Theorem:5.4} follows.
$\Box$
\begin{theorem}
\label{Theorem:5.5}
Let $p=p_{N+1}<2$. 
Let $\beta>0$ and set $\Phi(\tau):=\tau[\log(e+\tau)]^\beta$ for $\tau\ge 0$. 
For any $T>0$, there exists $\gamma=\gamma(N,\beta,T)>0$ such that, 
if $h\in{\mathcal L'}$ satisfies 
\begin{equation}
\label{eq:5.9}
 \sup_{x'\in{\bf R}^{N-1}}\int_{B'(x',\sigma)}\Phi(T^{\frac{1}{p-1}}h(y'))\,dy
 \le\gamma T^{\frac{N-1}{2}}\biggr[\log\left(e+\frac{\sqrt{T}}{\sigma}\right)\biggr]^{\beta-\frac{N+1}{2}}
\end{equation}
for all $\sigma\in(0,\sqrt{T})$, 
then problem~\eqref{eq:P} with $\mu=h(x')\otimes\delta_1(x_N)$
possesses a solution~$u$ in $Q_T$, with $u$ satisfying 
\[
0\le u(x,t)
\le C\frac{x_N}{t} \Gamma_1(x_N,t) 
\Psi^{-1}\left[\Gamma_{N-1}(t)\Psi(T^{\frac{1}{p-1}}h)\right](x)
\quad \mbox{for a.a.~$(x,t)\in Q_T$},
\]
for some $C>0$. 
\end{theorem}
{\bf Proof.}
Assume \eqref{eq:5.9}. We can assume, without loss of generality, that $\gamma\in(0,1)$. 
Define $\Psi$ as in the proof of Theorem \ref{Theorem:5.4}. 
Set 
\[
z(x,t) :=\left[\Gamma_{N-1}(t)\Psi(T^{\frac{1}{p-1}}h)\right](x).
\]
By Lemma~\ref{Lemma:2.2} and \eqref{eq:5.9} we see that 
\begin{equation*}
\begin{split}
\|z(t)\|_\infty
&\le C t^{-\frac{N-1}{2}} \sup_{x'\in{\bf R}^{N-1}} \int_{B'(x',\sqrt{t} )}\Psi(T^{\frac{1}{p-1}}h)\,dy'\\
&\le C\gamma t_T^{-\frac{N-1}{2}}|\log t_T|^{\beta-\frac{N+1}{2}}
\le C t_T^{-\frac{N-1}{2}}|\log t_T|^{\beta-\frac{N+1}{2}}
\quad\mbox{for $t\in(0,T)$},
\end{split}
\end{equation*}
where $t_T:=t/(2T)\in(0,1/2)$. 
By the same argument as in the proof of Theorem~\ref{Theorem:5.4} 
we have
\[
\begin{aligned}
	&0\leq A(z(x,t)) \leq 
	C\gamma^{p-1-\epsilon} t_T^{-\frac{N-1}{2}(p-1)}
	|\log t_T|^{(\beta-\frac{N+1}{2})(p-1)-\beta p}
	=
	C\gamma^{p-1-\epsilon} t_T^{-\frac{N-1}{N+1}} |\log t_T|^{-\beta-1}, \\
	&0\leq B(z(x,t)) \leq 
	C|\log t_T|^\beta, 
\end{aligned}
\]
for all $(x,t)\in Q_T$, where $C$ is independent of $\gamma$. It follows that 
\[
\begin{aligned}
	&\|B(z(t))\|_\infty
	\int_0^t s^{-(p-1)} \|A(z(s))\|_\infty\,ds  \\
	&= (2T)^{-(p-1)} \|B(z(t))\|_\infty
	\int_0^t s_T^{-(p-1)} \|A(z(s))\|_\infty\,ds \\
	&\leq 
	C T^{-(p-1)} \gamma^{p-1-\epsilon} |\log t_T|^\beta 
	\int_0^t s_T^{-1} |\log s_T|^{-\beta-1} \,ds 
	\leq 
	C T^{2-p} \gamma^{p-1-\epsilon}. 
\end{aligned}
\]
Then Theorem~\ref{Theorem:5.2} leads to  
the desired conclusion. The proof is complete.
$\Box$
\subsection{Optimal singularities}
Applying our necessary conditions and sufficient conditions on the solvability for problem~(P),
for any $z\in\overline{\Omega}$, we find a function $f_z\in{\mathcal L}$ 
with the following properties:
\begin{itemize}
  \item there exists $R>0$ such that $f_z$ is smooth in $B_\Omega(z,R)\setminus\{z\}$ and $f_z=0$ outside $B_\Omega(z,R)$;
  \item there exists $\kappa_z>0$ such that 
  problem~\eqref{eq:P} with $\mu=\kappa x_N f_z(x)$, where $\kappa>0$, possesses a local-in-time solution if $\kappa<\kappa_z$ 
  and it possesses no local-in-time solutions if $\kappa>\kappa_z$.
\end{itemize}
Similarly to Section~1.1, 
we term the singularity of the function $f_z$ at $x=z$ 
an {\it optimal singularity} of initial data for the solvability of problem~\eqref{eq:P} at $x=z$. 
\vspace{3pt}

We find optimal singularities of initial data for the solvability of problem~\eqref{eq:P} at $z\in\Omega$. 
By Corollary~\ref{Corollary:5.1}~(ii) it suffices to consider 
the case $p\ge p_N$. 
For any set $E$ in ${\mathbb R}^N$, we denote by $\chi_E$ the characteristic function of $E$. 
\begin{theorem}
\label{Theorem:5.6}
Let $z\in\Omega$.
Set 
\[
f_z(x):=
\left\{
\begin{array}{ll}
|x-z|^{-\frac{2}{p-1}}\chi_{B_\Omega(z,1)}(x) & \mbox{if}\quad p>p_N,\vspace{3pt}\\
|x-z|^{-N}\left|\log|x-z|\right|^{-\frac{N}{2}-1}\chi_{B_\Omega(z,1/2)}(x) 
& \mbox{if}\quad p=p_N,
\end{array}
\right.
\]
for $x\in\Omega$. 
Then there exists $\kappa_z>0$ with the following properties:
\begin{itemize}
  \item[{\rm (i)}] 
  problem~\eqref{eq:P} possesses a local-in-time solution with $\mu=\kappa x_Nf_z(x)$ if $0<\kappa<\kappa_z$;
  \item[{\rm (ii)}] 
  problem~\eqref{eq:P} possesses no local-in-time solutions with $\mu=\kappa x_Nf_z(x)$ if $\kappa>\kappa_z$.
\end{itemize}
Here $\displaystyle{\sup_{z\in\Omega}}\,\kappa_z<\infty$. 
\end{theorem}
{\bf Proof.}
Let $z=(z',z_N)\in\Omega$, $\kappa>0$, and $\mu=\kappa x_Nf_z(x)$ in ${\mathcal M}$. 
Assume that problem~\eqref{eq:P} possesses a local-in-time solution. 
By Theorem~\ref{Theorem:1.2} we have
\begin{equation}
\label{eq:5.10}
\kappa \int_{B(z,\sigma)}y_Nf_z(y)\,dy  
\leq C\sigma^{-\frac{2}{p-1}}\int_{B(z,\sigma)}y_N\,dy
\le Cz_N\sigma^{N-\frac{2}{p-1}}
\end{equation}
for all small enough $\sigma>0$. 
Furthermore, if $p=p_N$, then 
\begin{equation}
\label{eq:5.11}
\kappa\int_{B(z,\sigma)}y_Nf_z(y)\,dy\leq Cz_N|\log\sigma|^{-\frac{N}{2}}
\end{equation}
for all small enough $\sigma>0$. 
On the other hand, it follows that
$$
\int_{B(z,\sigma)}y_Nf_z(y)\,dy
\ge 
\left\{
\begin{array}{ll}
Cz_N\sigma^{N-\frac{2}{p-1}} & \mbox{if}\quad p>p_N,\vspace{3pt}\\
Cz_N|\log \sigma|^{-\frac{N}{2}} & \mbox{if}\quad p=p_N,
\end{array}
\right.
$$
for all small enough $\sigma>0$.
This together with \eqref{eq:5.10} and \eqref{eq:5.11} implies that $\kappa_z$ is uniformly bounded on $\Omega$. 

On the other hand, 
if $p>p_N$, then
we find $\alpha>1$ such that 
\[
\sup_{x\in\overline{\Omega}}\int_{B(x,\sigma)}\frac{y_N}{y_N+\sigma}(\kappa f_z(y))^\alpha\,dy
\le\kappa^\alpha\int_{B(z,\sigma)}|y-z|^{-\frac{2\alpha}{p-1}}\,dy\le C\kappa^\alpha\sigma^{N-\frac{2\alpha}{p-1}}
\]
for all $\sigma\in(0,1)$. 
If $p=p_N$, then, 
for any $\beta\in(0,N/2)$, we have
\begin{equation*}
\begin{split}
 & \sup_{x\in\overline{\Omega}}\int_{B(x,\sigma)}\kappa f_z(y)[\log(e+\kappa f_z(y))]^\beta\,dy\\
 & \le C\kappa\sup_{x\in\overline{\Omega}}\int_{B(z,\sigma)}|y-z|^{-N}|\log|y-z||^{-\frac{N}{2}-1+\beta}\,dy
\le C\kappa |\log\sigma|^{-\frac{N}{2}+\beta}
\end{split}
\end{equation*}
for all small enough $\sigma>0$ and $\kappa\in(0,1)$.
Then, if $\kappa>0$ is small enough, 
by Theorem~\ref{Theorem:5.3} with $p>p_N$ and 
Theorem~\ref{Theorem:5.4} with $\ell=0$  
we find a local-in-time solution to problem~\eqref{eq:P}. 
Therefore, thanks to Lemma~\ref{Lemma:2.5}, 
we find the desired constant $\kappa_z$, and the proof is complete. 
$\Box$
\vspace{3pt}

Next, we consider the case $z\in\partial\Omega$. 
Then it suffices to consider the case $z=0\in\partial\Omega$. 
\begin{theorem}
\label{Theorem:5.7}
Set 
\begin{equation*}
f(x):=
\left\{
\begin{array}{ll}
|x|^{-\frac{2}{p-1}}\chi_{B_\Omega(0,1)}(x)  & \mbox{if}\quad p>p_{N+1},\vspace{3pt}\\
|x|^{-N-1}\left|\log|x|\right|^{-\frac{N+1}{2}-1}\chi_{B_\Omega(0,1/2)}(x) 
& \mbox{if}\quad p=p_{N+1},
\end{array}
\right.
\end{equation*}
for $x\in\Omega$. 
Then there exists $\kappa_0>0$ with the following properties:
\begin{itemize}
  \item[{\rm (i)}] 
  problem~\eqref{eq:P} possesses a local-in-time solution with $\mu=\kappa x_Nf(x)$ if $0<\kappa<\kappa_0$;
  \item[{\rm (ii)}] 
  problem~\eqref{eq:P} possesses no local-in-time solutions with $\mu=\kappa x_Nf(x)$ if $\kappa>\kappa_0$.
\end{itemize}
\end{theorem}
{\bf Proof.}
Let $\kappa>0$ and $\mu=\kappa x_Nf(x)$ in ${\mathcal M}$.
Assume that problem~\eqref{eq:P} possesses a local-in-time solution.  
By Theorem~\ref{Theorem:1.2} we have
\begin{equation}
\label{eq:5.12}
\kappa\int_{B_\Omega(0,\sigma)}y_Nf(y)\,dy
\le C\sigma^{-\frac{2}{p-1}}\int_{B_\Omega(0,\sigma)}y_N\,dy
\le C\sigma^{N+1-\frac{2}{p-1}}
\end{equation}
for all small $\sigma>0$. 
Furthermore, if $p=p_{N+1}$, then 
\begin{equation}
\label{eq:5.13}
\kappa\int_{B_\Omega(0,\sigma)}y_Nf(y)\,dy\le C|\log\sigma|^{-\frac{N+1}{2}}
\end{equation}
for all small $\sigma>0$. 
It follows that 
\begin{equation}
\label{eq:5.14}
\int_{B_\Omega(0,\sigma)}y_Nf(y)\,dy
\ge
\left\{
\begin{array}{ll}
C\sigma^{N+1-\frac{2}{p-1}}\quad&\mbox{if}\quad p>p_{N+1},\vspace{3pt}\\
C|\log\sigma|^{-\frac{N+1}{2}}\quad&\mbox{if}\quad p=p_{N+1},
\end{array}
\right.
\end{equation}
for small enough $\sigma>0$. 
By \eqref{eq:5.12}, \eqref{eq:5.13}, and \eqref{eq:5.14} 
we see that $\kappa_0\le C$.

On the other hand, 
if $p>p_{N+1}$, 
we find $\alpha>1$ such that 
$$
\int_{B_\Omega(x,\sigma)}\frac{y_N}{y_N+\sigma}(\kappa f(y))^\alpha\,dy
\le \kappa^\alpha \sigma^{-1} \int_{B_\Omega(0,3\sigma)}y_N|y|^{-\frac{2\alpha}{p-1}}\,dy
\le C\kappa^\alpha \sigma^{N-\frac{2\alpha}{p-1}}
$$
for $x\in B_\Omega(0,2\sigma)$. 
Furthermore, we have
$$
\int_{B_\Omega(x,\sigma)}\frac{y_N}{y_N+\sigma}(\kappa f(y))^\alpha\,dy
\le \kappa^\alpha\int_{B_\Omega(x,\sigma)} |y|^{-\frac{2\alpha}{p-1}}\,dy
\le C\kappa^\alpha \sigma^N|x|^{-\frac{2\alpha}{p-1}}
\le C\kappa^\alpha \sigma^{N-\frac{2\alpha}{p-1}}
$$
for $x\in\overline{\Omega}\setminus B_\Omega(0,2\sigma)$. 
These imply that 
$$ 
\sup_{x\in\overline{\Omega}}\int_{B_\Omega(x,\sigma)}\frac{y_N}{y_N+\sigma}(\kappa f(y))^\alpha\,dy
\le C\kappa^\alpha\sigma^{N-\frac{2\alpha}{p-1}}
$$
for all $\sigma\in(0,1)$ if $p>p_{N+1}$.
If $p=p_{N+1}$, 
for any $\beta\in(0,N/2)$, we have
\begin{equation*}
\begin{split}
 & \sup_{x\in\overline{\Omega}}\int_{B(x,\sigma)}\kappa y_Nf(y)[\log(e+\kappa f(y))]^\beta\,dy\\
 & \le C\kappa\int_{B(0,\sigma)}|y|^{-N}|\log|y||^{-\frac{N}{2}-1+\beta}\,dy\le C\kappa|\log\sigma|^{-\frac{N}{2}+\beta}
\end{split}
\end{equation*}
for all small enough $\sigma>0$ and all $\kappa\in(0,1)$. 
Then, if $\kappa>0$ is small enough, 
by Theorem~\ref{Theorem:5.3} with $p>p_{N+1}$ and 
Theorem~\ref{Theorem:5.4} with $\ell=1$ 
we find a local-in-time solution to problem~\eqref{eq:P}. 
Therefore, thanks to Lemma~\ref{Lemma:2.5}, 
we find the desired constant $\kappa_0$, and the proof is complete. 
$\Box$
\vspace{5pt}

Similarly, we have:
\begin{theorem}
\label{Theorem:5.8}
Let $N\ge 2$ and $p_{N+1}\le p<2$. 
Set 
\[
h(x'):=
\left\{
\begin{array}{ll}
|x'|^{2-\frac{2}{p-1}}\chi_{B'(0,1)}(x') & \mbox{if}\quad p>p_{N+1},\vspace{3pt}\\
|x'|^{-N+1}\left|\log|x'|\right|^{-\frac{N+1}{2}-1}\chi_{B'(0,1/2)}(x') & \mbox{if}\quad p=p_{N+1},
\end{array}
\right.
\]
for $x'\in{\bf R}^{N-1}$. 
Consider problem~\eqref{eq:P} with 
$\mu=\kappa h(x')\otimes\delta_1(x_N)$, 
where $\kappa>0$. 
Then there exists $\kappa_0>0$ with the following properties:
\begin{itemize}
  \item[{\rm (i)}] 
  problem~\eqref{eq:P} possesses a local-in-time solution if $0<\kappa<\kappa_0$;
  \item[{\rm (ii)}] 
  problem~\eqref{eq:P} possesses no local-in-time solutions if $\kappa>\kappa_0$.
\end{itemize} 
\end{theorem}
{\bf Proof.}
Let $\kappa>0$ and $\mu=\kappa h(x')\otimes\delta_1(x_N)\in{\mathcal M}$. 
Assume that problem~\eqref{eq:P} possesses a local-in-time solution. 
By Theorem~\ref{Theorem:1.2} we have
\begin{align*}
C^{-1}\kappa\sigma^{N+1-\frac{2}{p-1}}
 & \le \kappa\int_{B'(0,\sigma)}h(y')\,dy'\le C\gamma\sigma^{N+1-\frac{2}{p-1}}\quad\mbox{if}\quad p>p_{N+1},\\
C^{-1}\kappa|\log\sigma|^{-\frac{N+1}{2}}
 & \le \kappa\int_{B'(0,\sigma)}h(y')\,dy'\le C|\log\sigma|^{-\frac{N+1}{2}}\quad\mbox{if}\quad p=p_{N+1},
\end{align*}
for all small enough $\sigma>0$. 
Then we see that $\kappa_0\le C$.

On the other hand, 
if $p>p_{N+1}$, 
we find $\alpha>1$ such that 
\begin{equation*}
\sup_{x'\in{\bf R}^{N-1}}\int_{B'(x',\sigma)}h(y)^\alpha\,dy
\le\kappa^\alpha\int_{B'(0,\sigma)}|y'|^{-2\alpha\frac{2-p}{p-1}}\,dy'
 \le C\kappa^\alpha\sigma^{N-1-2\alpha\frac{2-p}{p-1}}
\end{equation*}
for all $\sigma\in(0,1)$. 
If $p=p_{N+1}$, 
for any $\beta\in(0,N/2)$, we have
\begin{equation*}
\begin{split}
 & \sup_{x'\in{\bf R}^{N-1}}\int_{B'(x',\sigma)}\kappa h(y')[\log(e+\kappa h(y'))]^\beta\,dy'\\
 & \le C\kappa\int_{B'(0,\sigma)}|y'|^{-N+1}|\log|y'||^{-\frac{N}{2}-1+\beta}\,dy'
 \le C|\log\sigma|^{-\frac{N}{2}+\beta}
\end{split}
\end{equation*}
for all small enough $\sigma>0$ and all $\kappa\in(0,1)$. 
Then, if $\kappa>0$ is small enough, 
by Theorem~\ref{Theorem:5.3} with $p>p_{N+1}$ and 
Theorem~\ref{Theorem:5.5} we find a local-in-time solution to problem~\eqref{eq:P}.
Therefore, thanks to Lemma~\ref{Lemma:2.5}, 
we find the desired constant $\kappa_0$, and the proof is complete. 
$\Box$
\medskip

\noindent
{\bf Acknowledgment.}
The first and second authors were supported 
in part by JSPS KAKENHI Grant Number JP19H05599.
The third author was supported in part by JSPS KAKENHI Grant Numbers 
JP19K14567 and JP22H01131.

\end{document}